\documentclass[10pt,reqno,oneside]{amsart}
\usepackage{verbatim,amsmath,amssymb,cite,epsfig,subfigure,color}

\numberwithin{equation}{section}

\newcommand{\Real}{\mathbb R}
\newcommand{\N}{\mathbb N}
\newcommand{\Z}{\mathbb Z}

\newcommand{\E}{{\mathcal E}}
\newcommand{\C}{{\mathcal C}}
\newcommand{\U}{{\mathcal U}}
\newcommand{\Y}{{\mathcal Y}}

\newtheorem{theorem}{Theorem}[section]

\newtheorem{remark}{Remark}[section]

\usepackage{pdfsync}

\begin{document}

\title[Quasi-nonlocal approximation of linear and circular chains]
{Analysis of the quasi-nonlocal approximation of linear and circular chains in the plane}

\author{Pavel B\v{e}l\'{\i}k}
\author{Mitchell Luskin}

\address{Pavel B\v{e}l\'{\i}k\\
Mathematics Department\\
Augsburg College\\
2211 Riverside Avenue\\
Minneapolis, MN 55454\\
U.S.A.} \email{belik@augsburg.edu}

\address{Mitchell Luskin \\
School of Mathematics \\
University of Minnesota \\
206 Church Street SE \\
Minneapolis, MN 55455 \\
U.S.A.} \email{luskin@umn.edu}

\thanks{
This work was supported in part by DMS-0757355, DMS-0811039, the Institute for Mathematics and Its Applications, and the University of Minnesota Supercomputing Institute. This work was also supported by the Department of Energy under Award Number DE-SC0002085.
}

\date{\today}

\begin{abstract}
We give an analysis of the stability and displacement error for
linear and circular atomistic chains in the plane when the
atomistic energy is approximated by the Cauchy--Born continuum
energy and by the quasi-nonlocal atomistic-to-continuum
coupling energy. We consider atomistic energies that include
Lennard-Jones type nearest neighbor and next nearest neighbor
pair-potential interactions.

Previous analyses for linear chains have shown that the
Cauchy--Born and quasi-nonlocal approximations reproduce (up to
the order of the lattice spacing) the atomistic lattice
stability for perturbations that are constrained to the line of
the chain. However, we show that the Cauchy--Born and
quasi-nonlocal approximations give a finite increase for the
lattice stability of a linear or circular chain under
compression when general perturbations in the plane are
allowed. We also analyze the increase of the lattice stability
under compression when pair-potential energies are augmented by
bond-angle energies. Our estimates of the largest strain for
lattice stability (the critical strain) are sharp (exact up to
the order of the lattice scale).

We then use these stability estimates and  modeling error
estimates for the linearized Cauchy--Born and quasi-nonlocal
energies to give an optimal order (in the lattice scale) {\em a
priori} error analysis for the approximation of the atomistic
strain in $\ell^2_\varepsilon$ due to an external force.
\end{abstract}

\maketitle{\allowdisplaybreaks\thispagestyle{empty}

\section{Introduction}
The quasicontinuum (QC) method ~\cite{Ortiz:1995a} is a general
approach for coupling atomistic models with Cauchy--Born
continuum models to achieve both accuracy and efficiency. Many
authors have improved, extended, and analyzed the QC method and
related atomistic-to-continuum coupling
methods~\cite{makridakis10,Legoll:2005,brian10,gaviniorbital,LinP:2003a,badia:onAtCcouplingbyblending,LinP:2006a,Ortner:2008a,E:2006,legoll09,Dobson:2008a,Gunzburger:2008a,Gunzburger:2008b,Miller:2008,PrudhommeBaumanOden:2006}.

In this paper, we give a linearized analysis of the stability
and strain error in $\ell^2_\varepsilon$ for linear and
circular atomistic chains in the plane when the atomistic
energy is approximated by the Cauchy--Born continuum energy and
by the quasi-nonlocal atomistic-to-continuum coupling
energy~\cite{Shimokawa:2004}. We consider first atomistic
energies that include only Lennard-Jones type nearest neighbor
and next nearest neighbor pair-potential interactions, and we
then consider atomistic energies that also include bond-angle
interactions.

We chose the quasi-nonlocal atomistic-to-continuum coupling
energy because uniformly spaced linear chains (one-dimensional
lattices) are equilibria for the coupling energy just as they
are for the atomistic and Cauchy--Born continuum
energies~\cite{ortner:qnl1d,doblusort:qce.stab}. This property
is called patch test consistency. Patch test consistent
extensions of the quasi-nonlocal energy to finite range
interactions have been given
in~\cite{E:2006,shapeev,finite.range}.

A uniformly strained one-dimensional lattice modeled by a
Lennard-Jones type atomistic interaction loses stability when
the strain reaches a critical value (the critical strain). We
seek to estimate the critical strain for quasicontinuum
energies and to then compare them with the critical strain for
the atomistic energy. We define such estimates to be sharp if
they are exact up to the order of the lattice spacing. Sharp
lattice stability and optimal order (in the lattice scale)
strain error analyses of the one-dimensional quasi-nonlocal
approximation have been given
in~\cite{Dobson:2008b,ortner:qnl1d,doblusort:qce.stab,mingyang}.
In this paper, we give a sharp stability and optimal order
strain error analysis of the linearized problem for the
quasi-nonlocal approximation of some simple objective
structures that are generated by a single affine
mapping~\cite{james_06}. Such objective structures include
linear chains, circular chains, and helical chains. We focus on
linear and circular chains in the plane, and we note that our
analysis allows general planar perturbations.

Previous analyses of linear chains have shown that the quasi-nonlocal approximation reproduces the lattice stability for perturbations that are constrained to the line of the chain. However, chains can undergo buckling instabilities under compression when general planar perturbations are allowed. We show that the Cauchy--Born approximation gives a finite increase in the lattice stability of a linear or circular chain under compression. We also analyze the increase of the lattice stability under compression when pair-potential interactions are augmented by bond-angle energies.

We restrict our analysis to the classical quasi-nonlocal approximation of chains with next nearest neighbor interactions~\cite{Shimokawa:2004}. The analysis in~\cite{finite.range} can likely be utilized with the analysis in this paper to obtain a quasi-nonlocal analysis of linear and circular chains for finite range interactions, but the arguments would be considerably more technical than the analysis presented in this paper. We give our error analysis for the linearization about linear and circular chains. The linear analysis in this paper can also likely be extended to a nonlinear analysis by utilizing the inverse function theorem techniques developed in~\cite{ortner:qnl1d}, but the details of this analysis would also greatly increase the complexity of the analysis that we present. For simplicity and clarity of exposition, we will present our results and analysis for chains in $\Real^2$, but we note that our analysis generalizes directly to helical chains in $\Real^3$.

In Section~\ref{sec:definitions}, we define the energy of a chain of atoms and its Cauchy--Born approximation. In Section~\ref{sec:variations}, detailed derivations of the first and second variations of the atomistic and Cauchy--Born energies are presented. We give sharp lattice stability results for unconstrained periodic chains in Section~\ref{sec:stab_linear} and for circular chains in Section~\ref{sec:stab_circ}. We introduce a model for bond-angle energy in Section~\ref{sec:regularization} and derive results for its contribution to lattice stability.

In Section~\ref{sec:truncation}, we give estimates for the modeling error due to the Cauchy--Born approximation, and we then give an error analysis of the linearized problems in Sections~\ref{sec:error} and \ref{sec:error_bending} based on our stability and modeling error analyses. Finally, in Section~\ref{sec:qnl}, we define the quasi-nonlocal approximation for periodic chains and give a sharp stability and error analysis for its approximation of the atomistic model.

We summarize the results in Section \ref{sec:conclusion}.

\section{Definitions}
\label{sec:definitions}
We will consider chains of atoms $y=\{y_\ell\}_{\ell\in\Z}$, where $y_\ell\in\Real^2$ denotes the position of the $\ell$-th atom in a plane. We will assume that in the reference configuration the chain is straight and that the distance between neighboring atoms is $\varepsilon>0$. We will only focus on \emph{periodically repeating} chains with period $N\in\N$ such that $N\varepsilon=1$, where by periodically repeating we mean that $y_{\ell+N}=y_\ell+(y_N-y_0)$ for all $\ell\in\Z$. This definition allows us to treat \emph{closed} chains of $N$ atoms (for which $y_N=y_0$), and also chains that are not closed, but such that the shapes of the overall configurations repeat every $N$ atoms (such as straight chains). In 3-D, periodically repeating chains would include, for example, helical chains.

Thus, we can define the space of $N$-periodically repeating chains, $\mathcal{Y}$, by
\begin{equation*}
  \mathcal{Y}
  =
  \{y=\{y_\ell\}_{\ell\in\Z}:\ y_\ell\in\Real^2\text{ and }y_{\ell+N}=y_\ell+(y_N-y_0)\text{ for all }\ell\in\Z\}.
\end{equation*}
We will also use the space of $N$-periodic mean-zero displacements, $\mathcal{U}\subset\Y$, defined as
\begin{equation*}
  \mathcal{U}
  =
  \{u=\{u_\ell\}_{\ell\in\Z}:\ u_\ell\in\Real^2\text{ and }u_{\ell+N}=u_\ell\text{ for all }\ell\in\Z,\text{ and }\sum_{\ell=1}^Nu_\ell=0\}.
\end{equation*}
We will refer to such displacements as 2-D displacements.

Finally, we will also consider the subspaces $\tilde{\Y}\subset\Y$ and $\tilde{\U}\subset\U$ of one-dimensional periodic chains and displacements, respectively, defined as
\begin{equation*}
  \tilde{\Y}
  =
  \{y=\{y_\ell\}_{\ell\in\Z}\in\Y:\ y_\ell\cdot(0,1)=0\text{ for all }\ell\in\Z\}
\end{equation*}
and
\begin{equation*}
  \tilde{\U}
  =
  \{u=\{u_\ell\}_{\ell\in\Z}\in\U:\ u_\ell\cdot(0,1)=0\text{ for all }\ell\in\Z\}
\end{equation*}
when studying linear chains of atoms constrained so that the atoms can only move in the direction of the chain.

For $y\in\Y$, we define the \emph{backward} finite differences
\begin{equation*}
  y'_{\ell}
  =
  \frac{y_{\ell}-y_{\ell-1}}{\varepsilon},
  \quad
  y''_{\ell}
  =
  \frac{y'_{\ell}-y'_{\ell-1}}{\varepsilon}
  =
  \frac{y'_{\ell}-2y'_{\ell-1}+y'_{\ell-2}}{\varepsilon^2},
  \quad\dots,\quad
  y^{(n)}_\ell
  =
  \frac{y^{(n-1)}_\ell-y^{(n-1)}_{\ell-1}}{\varepsilon}
  \quad\text{ for }n\ge2,
\end{equation*}
and write $y'=\{y'_\ell\}_{\ell\in\Z}$, $y''=\{y''_\ell\}_{\ell\in\Z}$, etc.

For $v,w\in\Real^2$, we will write $v\cdot w$ for the usual inner product in $\Real^2$ and $\|v\|=\sqrt{v.v}$ for the usual Euclidean norm. For $v,w\in\Y$, we define the inner product
\begin{equation*}
  \langle v,w\rangle
  =
  \varepsilon\sum_{\ell=1}^Nv_\ell\cdot w_\ell,
\end{equation*}
and the associated norm
\begin{equation*}
  \|v\|_{\ell^2_\varepsilon}
  =
  \left(\varepsilon\sum_{\ell=1}^Nv_\ell\cdot v_\ell\right)^{1/2}
  =
  \left(\varepsilon\sum_{\ell=1}^N\|v_\ell\|^2\right)^{1/2}.
\end{equation*}
When estimating approximation errors, we will use the \emph{negative-order} norm
\begin{equation*}
  \|v\|_*
  =
  \sup_{w\in\U\setminus\{0\}}
  \frac{\langle v,w\rangle}{\|w'\|_{\ell^2_\varepsilon}}.
\end{equation*}

For a chain $y\in\Y$, we can now define the ``atomistic'' energy, $\E^{\text{a}}(y)$, via nearest neighbor and next nearest neighbor interactions
\begin{align}
  \label{eq:energy_atomistic}
  \E^{\text{a}}(y)
  &=
  \varepsilon\sum_{\ell=1}^N
  \left[
  \phi\left(\frac{\|y_{\ell}-y_{\ell-1}\|}{\varepsilon}\right)
  +
  \phi\left(\frac{\|y_{\ell+1}-y_{\ell-1}\|}{\varepsilon}\right)
  \right]\notag\\
  &=
  \varepsilon\sum_{\ell=1}^N
  \left[
  \phi\left(\|y'_{\ell}\|\right)
  +
  \phi\left(\|y'_{\ell+1}+y'_{\ell}\|\right)
  \right],
\end{align}
and its local Cauchy--Born approximation~\cite{Dobson:2008b,doblusort:qce.stab,yang.e.cb,arroyo04}, $\E^{\text{CB}}(y)$,
\begin{equation*}
  \E^{\text{CB}}(y)
  =
  \varepsilon\sum_{\ell=1}^N\E^{\text{CB}}_{\ell}(y),
\end{equation*}
where
\begin{equation*}
  \E^{\text{CB}}_{\ell}(y)
  =
  \frac{1}{2}\phi\left(\frac{\|y_{\ell+1}-y_{\ell}\|}{\varepsilon}\right)
  +
  \frac{1}{2}\phi\left(\frac{\|y_{\ell}-y_{\ell-1}\|}{\varepsilon}\right)
  +
  \frac{1}{2}\phi\left(2\frac{\|y_{\ell+1}-y_{\ell}\|}{\varepsilon}\right)
  +
  \frac{1}{2}\phi\left(2\frac{\|y_{\ell}-y_{\ell-1}\|}{\varepsilon}\right).
\end{equation*}
In the above approximation, we took into account only the nearest neighbor on either side of an atom at the position $y_\ell$ and extrapolated linearly to approximate the position of the next nearest neighbor on either side. Thus, the distance between this approximating next nearest neighbor and the original atom is twice the distance between the nearest neighbor and the original atom. Note that due to the chain being periodically repeating, we can write
\begin{equation}
  \label{eq:energy_cb}
  \E^{\text{CB}}(y)
  =
  \varepsilon\sum_{\ell=1}^N
  \left[
  \phi\left(\|y'_{\ell}\|\right)
  +
  \phi\left(2\|y'_{\ell}\|\right)
  \right].
\end{equation}

\section{Variations}
\label{sec:variations}
In the following sections, we will study the linearizations of the atomistic energy \eqref{eq:energy_atomistic} and Cauchy--Born energy \eqref{eq:energy_cb} about periodically repeating deformations $y\in\Y$ satisfying the particular constraints outlined below. We begin by deriving expressions for their first and second variations. We first define the strains
\begin{equation}
  \label{eq:strains}
  F_1
  =
  \frac{\|y_{\ell}-y_{\ell-1}\|}{\varepsilon}
  =
  \|y'_{\ell}\|
  \qquad\text{ and }\qquad
  F_2
  =
  \frac{\|y_{\ell+1}-y_{\ell-1}\|}{2\varepsilon}
  =
  \frac{\|y'_{\ell+1}+y'_{\ell}\|}{2},
\end{equation}
and we will assume that both $F_1$ and $F_2$ are independent of $\ell$. This assumption allows us to study, among others, linear (1-D) and circular deformations $y\in\Y$ with uniformly spaced atoms. Note that a 3-D helical chain also satisfies this assumption.

The first variation of the atomistic energy \eqref{eq:energy_atomistic} about a periodic configuration $y\in\Y$ satisfying \eqref{eq:strains} for all $\ell\in\Z$ is
\begin{align*}
  \delta\E^{\text{a}}(y)[u]
  &=
  \varepsilon\sum_{\ell=1}^N
  \left[
  \frac{\phi'(\|y'_{\ell}\|)}{\|y'_{\ell}\|}\,y'_{\ell}\cdot u'_{\ell}+\frac{\phi'(\|y'_{\ell+1}+y'_{\ell}\|)}{\|y'_{\ell+1}+y'_{\ell}\|}\,(y'_{\ell+1}+y'_{\ell})\cdot(u'_{\ell+1}+u'_{\ell})
  \right]\\
  &=
  \varepsilon\sum_{\ell=1}^N
  \left[
  \frac{\phi'(F_1)}{F_1}\,y'_{\ell}\cdot u'_{\ell}+\frac{\phi'(2F_2)}{2F_2}\,(y'_{\ell+1}+y'_{\ell})\cdot(u'_{\ell+1}+u'_{\ell})
  \right]\\
  &=
  \varepsilon\sum_{\ell=1}^N
  \left[
  \left(\frac{\phi'(F_1)}{F_1}+\frac{2\phi'(2F_2)}{F_2}\right)y'_{\ell}\cdot u'_{\ell}
  -
  \varepsilon^2\frac{\phi'(2F_2)}{2F_2}\,y''_{\ell}\cdot u''_{\ell}
  \right],
\end{align*}
where we have used the identity (with $A=I_2$, the $2\times2$ identity matrix)
\begin{equation}
  \label{eq:parallelogram}
  \left(y'_{\ell+1}+y'_{\ell}\right)\cdot A\left(u'_{\ell+1}+u'_{\ell}\right)
  =
  2\,y'_{\ell+1}\cdot Au'_{\ell+1}+2\,y'_{\ell}\cdot Au'_{\ell}-\varepsilon^2y''_{\ell+1}\cdot Au''_{\ell+1}
\end{equation}
and rearranged the sum using the periodicity of the chain.

The second variation of the atomistic energy \eqref{eq:energy_atomistic} is
\begin{align}
  \label{eq:E^a''_second}
  \delta^2\E^{\text{a}}(y)[u,v]
  &=
  \varepsilon\sum_{\ell=1}^N
  \left[
  \left(\phi''(\|y'_{\ell}\|)-\frac{\phi'(\|y'_{\ell}\|)}{\|y'_{\ell}\|}\right)\frac{u'_{\ell}\cdot y'_{\ell}}{\|y'_{\ell}\|}\frac{y'_{\ell}\cdot v'_{\ell}}{\|y'_{\ell}\|}\right.
  +
  \frac{\phi'(\|y'_{\ell}\|)}{\|y'_{\ell}\|}u'_{\ell}\cdot v'_{\ell}\notag\\
  &\quad\quad+
  \left(\phi''(\|y'_{\ell+1}+y'_{\ell}\|)-\frac{\phi'(\|y'_{\ell+1}+y'_{\ell}\|)}{\|y'_{\ell+1}+y'_{\ell}\|}\right)\frac{(u'_{\ell+1}+u'_{\ell})\cdot(y'_{\ell+1}+y'_{\ell})}{\|y'_{\ell+1}+y'_{\ell}\|}\frac{(y'_{\ell+1}+y'_{\ell})\cdot(v'_{\ell+1}+v'_{\ell})}{\|y'_{\ell+1}+y'_{\ell}\|}\notag\\
  &\quad\quad+
  \left.\frac{\phi'(\|y'_{\ell+1}+y'_{\ell}\|)}{\|y'_{\ell+1}+y'_{\ell}\|}(u'_{\ell+1}+u'_{\ell})\cdot(v'_{\ell+1}+v'_{\ell})\right]\notag\\
  &=
  \varepsilon\sum_{\ell=1}^N
  \left[
  u'_{\ell}\cdot\left(\phi''(F_1)P_\ell+\frac{\phi'(F_1)}{F_1}(I_2-P_\ell)\right)v'_{\ell}\right.\\
  &\quad\quad+
  \,(u'_{\ell+1}+u'_{\ell})\cdot\left(\phi''(2F_2)\tilde{P}_\ell+\frac{\phi'(2F_2)}{2F_2}(I_2-\tilde{P}_\ell)\right)(v'_{\ell+1}+v'_\ell)\biggr],\notag\\
  &=
  \label{eq:E^a''_first}
  \varepsilon\sum_{\ell=1}^N
  \left[
  u'_{\ell}\cdot\left(\phi''(F_1)P_\ell+\frac{\phi'(F_1)}{F_1}(I_2-P_\ell)\right)v'_{\ell}\right.\\
  &\quad\quad+
  \,(u'_{\ell+1}+u'_{\ell})\cdot\left(\phi''(2F_1)\tilde{P}_\ell+\frac{\phi'(2F_1)}{2F_1}(I_2-\tilde{P}_\ell)\right)(v'_{\ell+1}+v'_\ell)\notag\\
  &\quad\quad+
  \,(u'_{\ell+1}+u'_{\ell})\cdot\left((\phi''(2F_2)-\phi''(2F_1))\tilde{P}_\ell+\left(\frac{\phi'(2F_2)}{2F_2}-\frac{\phi'(2F_1)}{2F_1}\right)(I_2-\tilde{P}_\ell)\right)(v'_{\ell+1}+v'_\ell)\biggr],\notag
\end{align}
where we have used the projection operators
\begin{equation*}
 P_\ell
 =
 \frac{y'_{\ell}}{\|y'_{\ell}\|}\otimes\frac{y'_{\ell}}{\|y'_{\ell}\|}
 \quad\text{ and }\quad
 \tilde{P}_\ell
 =
 \dfrac{y'_{\ell+1}+y'_{\ell}}{\|y'_{\ell+1}+y'_{\ell}\|}\otimes\dfrac{y'_{\ell+1}+y'_{\ell}}{\|y'_{\ell+1}+y'_{\ell}\|}.
\end{equation*}
It is easy to see that $P^2_\ell=P_\ell$ and $\tilde{P}^2_\ell=\tilde{P}_\ell$ for all $\ell$. If we now apply identity \eqref{eq:parallelogram} to the term $(u'_{\ell+1}+u'_{\ell})\cdot\left(\phi''(2F_1)\tilde{P}_\ell\right)(v'_{\ell+1}+v'_\ell)$ in \eqref{eq:E^a''_first}, we obtain
\begin{align*}
  \delta^2\E^{\text{a}}(y)[u,v]
  &=
  \varepsilon\sum_{\ell=1}^N
  \left[
  u'_{\ell}\cdot\left(\phi''(F_1)P_\ell+\frac{\phi'(F_1)}{F_1}(I_2-P_\ell)\right)v'_{\ell}\right.\\
  &\qquad+
  (u'_{\ell+1}+u'_{\ell})\cdot\left(\frac{\phi'(2F_1)}{2F_1}(I_2-\tilde{P}_\ell)\right)(v'_{\ell+1}+v'_\ell)
  -
  \varepsilon^2u''_{\ell+1}\cdot\left(\phi''(2F_1)\tilde{P}_\ell\right)v''_{\ell+1}\\
  &\qquad+
  2\,u'_{\ell+1}\cdot\left(\phi''(2F_1)\tilde{P}_\ell\right)v'_{\ell+1}
  +
  2\,u'_{\ell}\cdot\left(\phi''(2F_1)\tilde{P}_\ell\right)v'_{\ell}\\
  &\qquad+
  \,(u'_{\ell+1}+u'_{\ell})\cdot\left((\phi''(2F_2)-\phi''(2F_1))\tilde{P}_\ell+\left(\frac{\phi'(2F_2)}{2F_2}-\frac{\phi'(2F_1)}{2F_1}\right)(I_2-\tilde{P}_\ell)\right)(v'_{\ell+1}+v'_\ell)\biggr].
\end{align*}
Finally, rearranging the sum, we obtain the final expression
\begin{align}
  \label{eq:E^a''}
  \delta^2\E^{\text{a}}(y)[u,v]
  &=
  \varepsilon\sum_{\ell=1}^N
  \left[
  u'_{\ell}\cdot\left(\bigr(\phi''(F_1)+4\,\phi''(2F_1)\bigr)P_\ell+\frac{\phi'(F_1)}{F_1}(I_2-P_\ell)\right)v'_{\ell}\right.\notag\\
  &\qquad+
  (u'_{\ell+1}+u'_{\ell})\cdot\left(\frac{\phi'(2F_1)}{2F_1}(I_2-\tilde{P}_\ell)\right)(v'_{\ell+1}+v'_\ell)
  -
  \varepsilon^2u''_{\ell+1}\cdot\left(\phi''(2F_1)\tilde{P}_\ell\right)v''_{\ell+1}\\
  &\qquad+
  2\,u'_{\ell}\cdot\left(\phi''(2F_1)(\tilde{P}_\ell+\tilde{P}_{\ell-1}-2P_\ell)\right)v'_{\ell}\notag\\
  &\qquad+
  \,(u'_{\ell+1}+u'_{\ell})\cdot\left((\phi''(2F_2)-\phi''(2F_1))\tilde{P}_\ell+\left(\frac{\phi'(2F_2)}{2F_2}-\frac{\phi'(2F_1)}{2F_1}\right)(I_2-\tilde{P}_\ell)\right)(v'_{\ell+1}+v'_\ell)\biggr].\notag
\end{align}

Similarly, the first variation of the Cauchy--Born approximation \eqref{eq:energy_cb} is
\begin{equation*}
  \delta\E^{\text{CB}}(y)[u]
  =
  \varepsilon\sum_{\ell=1}^N
  \frac{\phi'(F_1)+2\,\phi'(2F_1)}{F_1}\,y'_{\ell}\cdot u'_{\ell},
\end{equation*}
and the second variation is
\begin{equation}
  \label{eq:E^CB''}
  \delta^2\E^{\text{CB}}(y)[u,v]
  =
  \varepsilon\sum_{\ell=1}^N
  u'_{\ell}\cdot\left(\bigr(\phi''(F_1)+4\,\phi''(2F_1)\bigr)P_\ell+\frac{\phi'(F_1)+2\,\phi'(2F_1)}{F_1}\left(I_2-P_\ell\right)\right)v'_{\ell}.
\end{equation}

\section{Stability of a linear chain}
\label{sec:stab_linear}
In this section, we will consider a 1-D chain $y_F\in\tilde{\Y}$ of atoms with nearest neighbor interatomic spacing $F\varepsilon$ and study its stability with respect to two types of perturbations. First, we only consider displacements $u\in\tilde{\U}$ so that $y_F+u\in\tilde{\Y}$; this means that the atoms can only move in the direction of the chain and $y_F+u$ is still a 1-D chain. These results are given in Theorems \ref{thm:stab_CB_1D} and \ref{thm:stab_a_1D} and were first explicitly given in~\cite{doblusort:qce.stab}. In the second approach, we consider 2-D displacements that allow atoms to move out of the straight line. These results are given in Theorems \ref{thm:stab_CB_2D} and \ref{thm:stab_a_2D}. We then discuss and compare the results.

\subsection{Stability of a 1-D constrained chain}
Consider the 1-D atomic configuration $y_F\in\tilde{\Y}$ with interatomic spacing $F\varepsilon$. Note that in this case $F_1=F_2=F$ (see \eqref{eq:strains} for the definitions of $F_1$ and $F_2$). We then have the following stability results with respect to displacements that preserve the one-dimensionality of the chain (cf.~\cite{doblusort:qce.stab}).
\begin{theorem}
  \label{thm:stab_CB_1D}
  Let $y_F\in\tilde{\Y}$ denote the 1-D configuration of atoms with nearest neighbor interatomic spacing $F\varepsilon$. Then
  \begin{equation}
    \label{eq:stab_CB_1D1D_estimate}
    \inf_{u\in\mathcal{\tilde{U}}\setminus\{0\}}
    \frac{\delta^2\E^{\text{CB}}(y_F)[u,u]}{\|u'\|_{\ell^2_\varepsilon}^2}
    =
    \phi''(F)+4\,\phi''(2F).
  \end{equation}
\end{theorem}
\begin{proof}
  This follows immediately from \eqref{eq:E^CB''}, since for $u\in\mathcal{\tilde{U}}\setminus\{0\}$ we have $(I_2-P_\ell)u'_\ell=0$ for all $\ell$, and therefore $\delta^2\E^{\text{CB}}(y_F)[u,u]=\left(\phi''(F)+4\,\phi''(2F)\right)\|u'\|_{\ell^2_\varepsilon}^2$.
\end{proof}

\begin{theorem}
  \label{thm:stab_a_1D}
  Let $y_F\in\tilde{\Y}$ denote the 1-D configuration of atoms with nearest neighbor interatomic spacing $F\varepsilon$. If $\phi''(2F)\le0$, then
  \begin{equation*}
    \inf_{u\in\mathcal{\tilde{U}}\setminus\{0\}}
    \frac{\delta^2\E^{\text{a}}(y_F)[u,u]}{\|u'\|_{\ell^2_\varepsilon}^2}
    =
    \phi''(F)+4\,\phi''(2F)
    -
    \varepsilon^2\mu_\varepsilon\phi''(2F),
  \end{equation*}
  where
  \begin{equation*}
    \mu_\varepsilon
    =
    \inf_{u\in\mathcal{\tilde{U}}\setminus\{0\}}
    \frac{\|u''\|_{\ell^2_\varepsilon}}{\|u'\|_{\ell^2_\varepsilon}}
    =
    2\pi+\mathcal{O}(\varepsilon^2)
    \quad\text{ as }
    \varepsilon\to0.
  \end{equation*}
\end{theorem}
\begin{proof}
  This follows from \eqref{eq:E^a''}, since for $u\in\mathcal{\tilde{U}}\setminus\{0\}$ we have $(I_2-P_\ell)u'_\ell=(I_2-\tilde{P}_\ell)(u'_{\ell+1}+u'_\ell)=0$ and $\tilde{P}_\ell+\tilde{P}_{\ell-1}-2P_\ell=0$ for all $\ell$, and therefore
  \begin{equation*}
    \delta^2\E^{\text{a}}(y_F)[u,u]
    =
    \left(\phi''(F)+4\,\phi''(2F)\right)\|u'\|_{\ell^2_\varepsilon}^2
    -
    \varepsilon^2\phi''(2F)\|u''\|_{\ell^2_\varepsilon}^2.
  \end{equation*}
  The result follows by applying the identity \cite{SuliMayers,doblusort:qce.stab}
  \begin{equation*}
    \inf_{u\in\mathcal{\tilde{U}}\setminus\{0\}}
    \frac{\|u''\|_{\ell^2_\varepsilon}}{\|u'\|_{\ell^2_\varepsilon}}
    =
    \frac{2\sin{(\pi\varepsilon)}}{\varepsilon}
    =
    2\pi+\mathcal{O}(\varepsilon^2)
    \quad\text{ as }
    \varepsilon\to0.
  \end{equation*}
\end{proof}

\begin{remark}
  We recall that stability means that the infima in Theorems \ref{thm:stab_CB_1D} and \ref{thm:stab_a_1D} are positive. For a typical potential $\phi$, such as the Lennard-Jones potential, the requirement of positiveness provides an upper bound on the strain $F$ for which $y_F$ is stable in the respective model. That is, it provides a bound on the amount of stretching the 1-D chain can undergo and remain stable. Note that the atomistic model exhibits slightly more stability than the Cauchy--Born approximation in the sense that the upper bound on $F$ for which $y_F$ is stable is larger in the atomistic model than in the Cauchy--Born model. Asymptotically, however, as $\varepsilon\to0$, the two infima agree. Finally, we also note that the stability region is only bounded by $0$ from below, that is, there is technically no bound on the amount of compression the 1-D chain can undergo and become unstable in either model.
\end{remark}

\subsection{Stability of a 1-D unconstrained chain}
Consider again the 1-D atomic configuration $y_F\in\tilde{\Y}$ with nearest neighbor interatomic spacing $F\varepsilon$ so that again $F_1=F_2=F$. Since under compression of the chain one would expect the chain to exhibit some type of buckling, we next consider displacements $u\in\U$ that allow the chain to become two-dimensional, and we provide stability results with respect to such displacements.
\begin{theorem}
  \label{thm:stab_CB_2D}
  Let $y_F\in\tilde{\Y}$ denote the 1-D configuration of atoms with nearest neighbor interatomic spacing $F\varepsilon$. Then
  \begin{equation*}
    \inf_{u\in\mathcal{U}\setminus\{0\}}
    \frac{\delta^2\E^{\text{CB}}(y_F)[u,u]}{\|u'\|_{\ell^2_\varepsilon}^2}
    =
    \min{\left\{\phi''(F)+4\,\phi''(2F),\,\dfrac{\phi'(F)+2\,\phi'(2F)}{F}\right\}}.
  \end{equation*}
\end{theorem}
\begin{proof}
  Using expression \eqref{eq:E^CB''} for the second variation of $\E^{\text{CB}}$, we immediately get
  \begin{align}
    \label{eq:stab_CB_1D2D_estimate}
    \delta^2\E^{\text{CB}}(y_F)[u,u]
    &=
    \varepsilon\sum_{\ell=1}^N\biggr[
    \left(\phi''(F)+4\,\phi''(2F)\right)\|P_\ell u'_\ell\|^2
    +
    \frac{\phi'(F)+2\,\phi'(2F)}{F}\|(I_2-P_\ell)u'_\ell\|^2\biggr]\notag\\
    &\ge
    \min{\left\{\phi''(F)+4\,\phi''(2F),\,\frac{\phi'(F)+2\,\phi'(2F)}{F}\right\}}
    \|u'\|_{\ell^2_\varepsilon}^2.
  \end{align}
  To obtain the expression for the infimum, we will consider two types of displacements. First, for any 1-D displacement $\tilde{u}\in\mathcal{\tilde{U}}$, we have $(I_2-P_\ell)\tilde{u}'_\ell=0$ and
  \begin{equation*}
    \delta^2\E^{\text{CB}}(y_F)[\tilde{u},\tilde{u}]
    =
    \left(\phi''(F)+4\,\phi''(2F)\right)
    \|\tilde{u}'\|_{\ell^2_\varepsilon}^2.
  \end{equation*}
  Next, consider the case of $N$ even and the displacement $\hat u\in{\U}$ such that for all $\ell$ we have $\hat{u}_{2\ell}=Cv$ and $\hat{u}_{2\ell+1}=-Cv$ for some $C>0$ and a vector $v\in\Real^2$ orthogonal to the chain $y_F$. This case corresponds to creating a zig-zag deformation of the 1-D chain. In this case it is easy to see that $P_\ell\hat{u}'_\ell=0$ for all $\ell$ and
  \begin{equation}
    \label{eq:variation_uhat_linear}
    \delta^2\E^{\text{CB}}(y_F)[\hat{u},\hat{u}]
    =
    \frac{\phi'(F)+2\,\phi'(2F)}{F}
    \|\hat{u}'\|_{\ell^2_\varepsilon}^2.
  \end{equation}
  Finally, in the case of odd $N$, the same zig-zag deformation $\hat{u}$ defined for $\ell=1,\dots,N-1$ and $\hat{u}_N=0$ will still satisfy $P_\ell\hat{u}'_\ell=0$ for all $\ell$ and \eqref{eq:variation_uhat_linear} still holds.
\end{proof}

For the stability in the atomistic model we have the following result.
\begin{theorem}
  \label{thm:stab_a_2D}
  Let $y_F\in\tilde{\Y}$ denote the 1-D configuration of atoms with nearest neighbor interatomic spacing $F\varepsilon$. If $\phi'(2F)\ge0$ and $\phi''(2F)\le0$, then, as $\varepsilon\to0$,
  \begin{equation*}
    \inf_{u\in\mathcal{U}\setminus\{0\}}
    \frac{\delta^2\E^{\text{a}}(y_F)[u,u]}{\|u'\|_{\ell^2_\varepsilon}^2}
    =
    \begin{cases}
      \min{\left\{\phi''(F)+4\,\phi''(2F)+\mathcal{O}(\varepsilon^2),\,\dfrac{\phi'(F)}{F}\right\}} & \text{ if $N$ is even},\\
      \min{\left\{\phi''(F)+4\,\phi''(2F)+\mathcal{O}(\varepsilon^2),\,\dfrac{\phi'(F)}{F}+\mathcal{O}(\varepsilon)\right\}} & \text{ if $N$ is odd}.\\
    \end{cases}
  \end{equation*}
\end{theorem}
\begin{proof}
  Using expression \eqref{eq:E^a''} for the second variation of $\E^{\text{a}}$ and observing that $\tilde{P}_\ell+\tilde{P}_{\ell-1}-2P_\ell=0$, we get
  \begin{align*}
    \delta^2\E^{\text{a}}(y_F)[u,u]
    &=
    \varepsilon\sum_{\ell=1}^N\biggr[
    \bigr(\phi''(F)+4\,\phi''(2F)\bigr)\|P_\ell u'_\ell\|^2
    +
    \frac{\phi'(F)}{F}\|(I_2-P_\ell)u'_\ell\|^2\\
    &\qquad\qquad+
    \frac{\phi'(2F)}{2F}\|(I_2-\tilde{P}_\ell)(u'_{\ell+1}+u'_\ell)\|^2
    -
    \varepsilon^2\phi''(2F)\|\tilde{P}_\ell u''_{\ell+1}\|^2\biggr]\\
    &\ge
    \min{\left\{\phi''(F)+4\,\phi''(2F),\,\dfrac{\phi'(F)}{F}\right\}}\|u'\|_{\ell^2_\varepsilon}^2.
  \end{align*}
  To obtain the expression for the infimum, one can again use the displacements $\tilde{u}$ and $\hat{u}$ as in the proof of Theorem \ref{thm:stab_CB_2D}. It is easy to verify that they satisfy (for all values of $\ell$ for $\tilde{u}$, for all values of $\ell$ for $\hat{u}$ if $N$ is even, and for all but three values of $\ell$ for $\hat{u}$ if $N$ is odd)
  \begin{gather*}
    (I_2-P_\ell)\tilde{u}'_\ell=(I_2-\tilde{P}_\ell)(\tilde{u}'_{\ell+1}+\tilde{u}'_{\ell-1})=0,
    \quad
    \tilde{P}_\ell\tilde{u}''_{\ell+1}=\tilde{u}''_{\ell+1},\\
    P_\ell\hat{u}'_\ell=\hat{u}'_{\ell+1}+\hat{u}'_{\ell}=\tilde{P}_\ell\hat{u}''_{\ell+1}=0.
  \end{gather*}
  In the case of odd $N$, we do not have $\hat{u}'_{\ell+1}+\hat{u}'_{\ell}=0$ for $\ell=1$, $N-1$, and $N$, only $\tilde{P}_\ell(\hat{u}'_{\ell+1}+\hat{u}'_{\ell})=0$, and therefore
  \begin{gather*}
    \delta^2\E^{\text{a}}(y_F)[\tilde{u},\tilde{u}]
    =
    \left(\phi''(F)+4\,\phi''(2F)\right)
    \|\tilde{u}'\|_{\ell^2_\varepsilon}^2
    -
    \left(\varepsilon^2\phi''(2F)\right)\varepsilon\sum_{\ell=1}^N\|\tilde{P}_\ell u''_{\ell+1}\|^2,\\
    \delta^2\E^{\text{a}}(y_F)[\hat{u},\hat{u}]
    =
    \frac{\phi'(F)}{F}\|\hat{u}'\|_{\ell^2_\varepsilon}^2
    \quad\text{ if $N$ is even,}\\
    \delta^2\E^{\text{a}}(y_F)[\hat{u},\hat{u}]
    =
    \frac{\phi'(F)}{F}\|\hat{u}'\|_{\ell^2_\varepsilon}^2
    +
    \mathcal{O}(\varepsilon)
    \quad\text{ if $N$ is odd,}\\
  \end{gather*}
  and the conclusion of the theorem follows.
\end{proof}

\begin{remark}
  \label{rem:angle}
  Comparing Theorems \ref{thm:stab_CB_2D} and \ref{thm:stab_a_2D} to their constrained one-dimensional counterparts \ref{thm:stab_CB_1D} and \ref{thm:stab_a_1D}, we see that in order for $y_F$ to be stable in the two-dimensional models, additional inequalities, $\dfrac{\phi'(F)+2\,\phi'(2F)}{F}>0$ and $\dfrac{\phi'(F)}{F}>0$, respectively, must now be satisfied in the limit as $\varepsilon\to0$. For typical potentials $\phi$, these inequalities provide lower bounds on the stretch $F$, below which the chain could undergo a zig-zag buckling as demonstrated in the proofs. Note that a zig-zag configuration with all nearest neighbor interatomic distances equal to $\varepsilon$ and all turning angles alternatingly equal to $\pm2\pi/3$ would produce a global minimum of the atomic energy $\E^{\text{a}}$. (Following~\cite{oneil}, we consider the signed turning angle $-\pi<\beta_\ell<\pi$ at each atom $y_\ell$ defined as the angle between $y'_\ell$ and $y'_{\ell+1}$, measured in the sense that negative sign of $\beta_\ell$ corresponds to a clockwise turn and positive sign to a counterclockwise turn.)

  We note that in all four theorems, the inequality determining the upper bound on $F$ is asymptotically the same, $\phi''(F)+4\,\phi''(2F)>0$. Intuitively, this is reasonable to expect, since under tension, the atoms would tend to align themselves along a straight line, thus erasing the lowest-order difference between the atomistic and Cauchy--Born models, and between the one-dimensional and two-dimensional models.

  On the other hand, there is a difference between the inequalities for the lower bound on $F$ in Theorems \ref{thm:stab_CB_2D} and \ref{thm:stab_a_2D}, $\dfrac{\phi'(F)+2\,\phi'(2F)}{F}>0$ and $\dfrac{\phi'(F)}{F}>0$, respectively.
  These lower bounds provide buckling thresholds for the one-dimensional chains, and it is reasonable to expect them to be different in the atomistic and the Cauchy--Born models, since the instability mode, the zig-zag deformation, cannot be described in the Cauchy--Born model.
  Since for typical potentials $\phi$ one has $\phi'(2F)>0$, the Cauchy--Born model exhibits more stability under compression compared to the atomistic model. The lower bound on $F$ in the Cauchy--Born model corresponds to $F<1$, for which also $\phi'(F)<0$. If one creates a small zig-zag perturbation in this configuration, the nearest neighbor distance increases, thus lowering the contribution to the energy from the nearest neighbor interaction. The next nearest neighbor distance also increases in the Cauchy--Born model, resulting in an increase in the overall energy. However, in the atomistic model, the next nearest neighbor distance stays the same, and the overall energy is thus decreased by the zig-zag perturbation. Thus the requirement of $\phi'(F)>0$ in the atomistic model seems quite reasonable.
\end{remark}

\section{Stability of a circular chain}
\label{sec:stab_circ}
In this section, we will consider uniform circular configurations $y_F\in\Y$ of $N$ atoms with nearest neighbor interatomic spacing $F\varepsilon$ and study their stability with respect to 2-D perturbations $u\in\U$. The atoms lie on a circle of radius $R$ that satisfies
\begin{equation*}
  F\varepsilon
  =
  2R\sin(\pi\varepsilon).
\end{equation*}
The distance between next nearest neighbors of the chain $y$ is $2F\varepsilon\cos(\pi\varepsilon)$. Therefore, we have
\begin{equation*}
  F_1=F
  \quad\text{ and }\quad
  F_2=F\cos\dfrac{\pi}{N}=F\cos{(\pi\varepsilon)}.
\end{equation*}

We note that due to the symmetry of the circle the only forces on the atoms in both the atomistic and the Cauchy--Born model are in the radial direction, i.e., the direction normal to the circle. It is straightforward to obtain from equations \eqref{eq:energy_atomistic} and \eqref{eq:energy_cb} that these forces vanish in the Cauchy--Born model if
\begin{equation*}
  \phi'(F)+2\,\phi'(2F)=0,
\end{equation*}
and they vanish in the atomistic model if
\begin{equation*}
  \phi'(F_1)+2\cos(\pi\varepsilon)\,\phi'(2F_2)=0
  \quad\text{ or }\quad
  \phi'(F)+2\cos(\pi\varepsilon)\,\phi'(2\cos(\pi\varepsilon)F)=0.
\end{equation*}
We remark that the above equations are special cases of the vanishing of the first variations of the respective energies, $\delta\E(y_F)\equiv 0$, as given in Section~\ref{sec:variations}. For Lennard-Jones type potentials $\phi$, these equations always have a solution, and thus both the atomistic and the Cauchy--Born model possess a circular equilibrium, although in general with different radii $R^a$ and $R^{\text{CB}}$. Using the inverse function theorem, the difference between the radii can be seen to be $R^a-R^{\text{CB}}=\mathcal{O}(\varepsilon^2)$ as $\varepsilon\to0$.

We now have the following result for the stability in the Cauchy--Born model.
\begin{theorem}
  \label{thm:stab_CB}
  Let $y_F\in\Y$ denote the uniform circular configuration of $N$ atoms with nearest neighbor interatomic spacing $F\varepsilon$. Then, as $\varepsilon\to0$,
  \begin{equation*}
    \inf_{u\in\mathcal{U}\setminus\{0\}}
    \frac{\delta^2\E^{\text{CB}}(y_F)[u,u]}{\|u'\|_{\ell^2_\varepsilon}^2}
    =
    \begin{cases}
      \min{\left\{\phi''(F)+4\,\phi''(2F),\,\dfrac{\phi'(F)+2\,\phi'(2F)}{F}\right\}} & \text{ if $N$ is even},\\
      \min{\left\{\phi''(F)+4\,\phi''(2F),\,\dfrac{\phi'(F)+2\,\phi'(2F)}{F}+\mathcal{O}(\varepsilon)\right\}} & \text{ if $N$ is odd}.\\
    \end{cases}
  \end{equation*}
\end{theorem}
\begin{proof}
  Using expression \eqref{eq:E^CB''} for the second variation of $\E^{\text{CB}}$, we immediately get
  \begin{align}
    \label{eq:stab_CB_2D_estimate}
    \delta^2\E^{\text{CB}}(y_F)[u,u]
    &=
    \varepsilon\sum_{\ell=1}^N\biggr[
    \left(\phi''(F)+4\,\phi''(2F)\right)\|P_\ell u'_\ell\|^2
    +
    \frac{\phi'(F)+2\,\phi'(2F)}{F}\|(I_2-P_\ell)u'_\ell)\|^2\biggr]\notag\\
    &\ge
    \min{\left\{\phi''(F)+4\,\phi''(2F),\,\frac{\phi'(F)+2\,\phi'(2F)}{F}\right\}}
    \|u'\|_{\ell^2_\varepsilon}^2.
  \end{align}
  To obtain the expression for the infimum, we will again consider two types of displacements. First, if $\tilde{u}_\ell=Cy_\ell$ for some $C>0$ and for all $\ell$, a case that corresponds to the pure expansion of the circle, then $\tilde{u}\in\U$, $(I_2-P_\ell)\tilde{u}'_\ell=0$, and
  \begin{equation*}
    \delta^2\E^{\text{CB}}(y_F)[\tilde{u},\tilde{u}]
    =
    \left(\phi''(F)+4\,\phi''(2F)\right)
    \|\tilde{u}'\|_{\ell^2_\varepsilon}^2.
  \end{equation*}
  Next, consider the case of $N$ even and the displacement $\hat{u}$ such that $\hat{u}_{2\ell}=Cy_{2\ell}$ and $\hat{u}_{2\ell+1}=-Cy_{2\ell+1}$ for some $C>0$ and for all $\ell$. This case corresponds to creating a zig-zag deformation of the circle. In this case it is easy to see that $\hat{u}\in\U$, $P_\ell\hat{u}'_\ell=0$, and
  \begin{equation}
    \label{eq:variation_uhat_circle}
    \delta^2\E^{\text{CB}}(y_F)[\hat{u},\hat{u}]
    =
    \frac{\phi'(F)+2\,\phi'(2F)}{F}
    \|\hat{u}'\|_{\ell^2_\varepsilon}^2.
  \end{equation}
  Finally, in the case of odd $N$, the same zig-zag deformation $\hat{u}$ defined for $\ell=1,\dots,N-1$ and $\hat{u}_N$ defined so that $\sum_{\ell=1}^N\hat{u}_\ell=0$ will have two segments, one joining the first and $N$-th atoms and one joining the $(N-1)$-st and $N$-th atoms, for which $P_\ell\hat{u}'_\ell\ne0$, thus creating a $\mathcal{O}(\varepsilon)$ perturbation in \eqref{eq:variation_uhat_circle}.
\end{proof}

Before we address the stability in the atomistic model, we note that for circular arrangements $\tilde{P}_\ell+\tilde{P}_{\ell-1}-2P_\ell\ne0$ for any $\ell$. However, it is easy to check that due to the geometry of the circle, we have
\begin{equation}
  \label{eq:projections_2}
  \|(\tilde{P}_\ell+\tilde{P}_{\ell-1}-2P_\ell)w\|=2\|w\|\sin^2\frac{\pi}{N}=2\|w\|\sin^2{(\pi\varepsilon)}\le2\pi^2\varepsilon^2\|w\|
  \quad\text{ for all }w\in\Real^2.
\end{equation}
We now have the following stability result.
\begin{theorem}
  \label{thm:stab_a}
  Let $y_F\in\Y$ denote the uniform circular configuration of $N$ atoms with nearest neighbor interatomic spacing $F\varepsilon$. If $\phi\in\mathcal{C}^3(0,\infty)$, $\phi'(2F)\ge0$, and $\phi''(2F)\le0$, then, as $\varepsilon\to0$,
  \begin{equation*}
    \inf_{u\in\mathcal{U}\setminus\{0\}}
    \frac{\delta^2\E^{\text{a}}(y_F)[u,u]}{\|u'\|_{\ell^2_\varepsilon}^2}
    =
    \begin{cases}
      \min{\left\{\phi''(F)+4\,\phi''(2F),\,\dfrac{\phi'(F)}{F}\right\}}+\mathcal{O}(\varepsilon^2) & \text{ if $N$ is even},\\
      \min{\left\{\phi''(F)+4\,\phi''(2F),\,\dfrac{\phi'(F)}{F}+\mathcal{O}(\varepsilon)\right\}}+\mathcal{O}(\varepsilon^2) & \text{ if $N$ is odd}.\\
    \end{cases}
  \end{equation*}
\end{theorem}
\begin{proof}
  First, recall the definitions $F_1=F$ and $F_2=F\cos(\pi\varepsilon)$, so $|F_2-F|\le\dfrac{\varepsilon^2}{2}$. Using the smoothness of $\phi$, there exists a constant $C_\phi>0$, independent of $\varepsilon$, such that
  \begin{equation}
    \label{eq:C_phi}
    \max\left\{\left|\phi''(2F_2)-\phi''(2F)\right|,\left|\frac{\phi'(2F_2)}{2F_2}-\frac{\phi'(2F)}{2F}\right|\right\}
    \le
    C_\phi\,\varepsilon^2.
  \end{equation}
  Using expression \eqref{eq:E^a''} for the second variation of $\E^{\text{a}}$, we now have
  \begin{align*}
    \delta^2\E^{\text{a}}(y_F)[u,u]
    &=
    \varepsilon\sum_{\ell=1}^N\biggr[
    \bigr(\phi''(F)+4\,\phi''(2F)\bigr)\|P_\ell u'_\ell\|^2
    +
    \frac{\phi'(F)}{F_1}\|(I_2-P_\ell)u'_\ell\|^2\\
    &\qquad\qquad+
    \frac{\phi'(2F)}{2F}\|(I_2-\tilde{P}_\ell)(u'_{\ell+1}+u'_\ell)\|^2
    -
    \varepsilon^2\phi''(2F)\|\tilde{P}_\ell u''_{\ell+1}\|^2\\
    &\qquad\qquad+
    2\,u'_{\ell}\cdot\left(\phi''(2F)(\tilde{P}_\ell+\tilde{P}_{\ell-1}-2P_\ell)\right)u'_{\ell}\\
    &\qquad\qquad+
    (\phi''(2F_2)-\phi''(2F)\|\tilde{P}_{\ell}(u'_{\ell+1}+u'_\ell)\|^2\\
    &\qquad\qquad+
    \left(\frac{\phi'(2F_2)}{2F_2}-\frac{\phi'(2F)}{2F}\right)\|(I_2-\tilde{P}_{\ell})(u'_{\ell+1}+u'_\ell)\|^2
    \biggr].
  \end{align*}
  Using the assumptions $\phi'(2F)\ge0$ and $\phi''(2F)\le0$, applying \eqref{eq:projections_2}, and \eqref{eq:C_phi} together with the triangle inequality, we get
  \begin{align}
    \label{eq:stab_a_2D_estimate}
    \delta^2\E^{\text{a}}(y_F)[u,u]
    &\ge
    \min{\left\{\phi''(F)+4\,\phi''(2F),\frac{\phi'(F)}{F}\right\}}\|u'\|_{\ell^2_\varepsilon}^2\notag\\
    &\quad-
    4\,\pi^2\varepsilon^2|\phi''(2F)|\,\|u'\|_{\ell^2_\varepsilon}^2
    -
    4\,\varepsilon^2C_\phi\|u'\|_{\ell^2_\varepsilon}^2.
  \end{align}
  To obtain the expression for the infimum, one can again use the displacements $\tilde{u}$ and $\hat{u}$ as in the proof of Theorem \ref{thm:stab_CB}. It is easy to verify that they satisfy
  \begin{gather*}
   \tilde{P}_\ell\tilde{u}''_{\ell+1}=0,
   \quad
   (I_2-\tilde{P}_\ell)(\tilde{u}'_{\ell+1}+\tilde{u}'_{\ell-1})=0,\\
   \tilde{P}_\ell\hat{u}''_{\ell+1}=0,
   \quad
   (I_2-\tilde{P}_\ell)(\hat{u}'_{\ell+1}+\hat{u}'_{\ell-1})=0,
  \end{gather*}
  except for two values of $\hat{u}_\ell$ when $N$ is odd, and the conclusion of the theorem follows.
\end{proof}

\begin{remark}
  Asymptotically, as $\varepsilon\to0$, the results in Theorems \ref{thm:stab_CB} and \ref{thm:stab_a} are equivalent to those in Theorems \ref{thm:stab_CB_2D} and \ref{thm:stab_a_2D}. The requirement of stability gives two inequalities, and these inequalities produce an upper and lower bound of $F$ for typical potentials $\phi$. Yet again, we can interpret the upper bound as a limit on the amount of stretching the chain can undergo before fracturing, while the lower bound can be interpreted as a limit on the amount of compression before the chain starts buckling and creating zig-zag segments.
\end{remark}

\section{Stabilization by a bond-angle energy}
\label{sec:regularization}
Many-body empirical potentials generally include the effect of bond angle in addition to two-body interactions~\cite{baskes07,tersoff88}.  We will study a simple bond-angle energy in this section and show that it suppresses the buckling modes under compression when added to the atomistic energy \eqref{eq:energy_atomistic} and the Cauchy--Born energy \eqref{eq:energy_cb}.

We will consider a bond-angle energy of the form
\begin{equation}
  \label{eq:E^b}
  \E^{\text{b}}(y)
  =
  \varepsilon\sum_{\ell=1}^N
  \alpha(1-\cos{\beta_\ell}),
\end{equation}
where $\alpha>0$ is a constant and $-\pi<\beta_\ell<\pi$ is the signed turning angle of the chain $y$ at the atom $y_\ell$. Recall from Remark~\ref{rem:angle} that $\beta_\ell$ is the angle between $y'_{\ell}$ and $y'_{\ell+1}$ and that negative sign of $\beta_\ell$ corresponds to a clockwise turn and positive sign to a counterclockwise turn. We note that
\begin{equation*}
  \cos{\beta_\ell}
  =
  \frac{y'_{\ell+1}}{\|y'_{\ell+1}\|}\cdot\frac{y'_\ell}{\|y'_\ell\|}.
\end{equation*}
To see the effect of this term on the overall energy, we first compute the first and second variations of $\E^{\text{b}}$.

For the first variation, we have the following result.
\begin{align}
  \label{eq:variation_1_Eb}
  \delta\E^{\text{b}}(y)[u]
  &=
  \alpha\varepsilon\sum_{\ell=1}^N
  \biggr[\left(-\frac{u'_{\ell+1}}{\|y'_{\ell+1}\|}+\left(\frac{u'_{\ell+1}}{\|y'_{\ell+1}\|}\cdot\frac{y'_{\ell+1}}{\|y'_{\ell+1}\|}\right)\frac{y'_{\ell+1}}{\|y'_{\ell+1}\|}\right)\cdot\frac{y'_\ell}{\|y'_\ell\|}\notag\\
  &\qquad\qquad+
  \left(-\frac{u'_\ell}{\|y'_\ell\|}+\left(\frac{u'_\ell}{\|y'_\ell\|}\cdot\frac{y'_\ell}{\|y'_\ell\|}\right)\frac{y'_\ell}{\|y'_\ell\|}\right)\cdot\frac{y'_{\ell+1}}{\|y'_{\ell+1}\|}\biggr]\\
  &=
  \alpha\varepsilon\sum_{\ell=1}^N
  \biggr[\cos{\beta_\ell}\left(\frac{u'_{\ell+1}}{\|y'_{\ell+1}\|}\cdot\frac{y'_{\ell+1}}{\|y'_{\ell+1}\|}+\frac{u'_\ell}{\|y'_\ell\|}\cdot\frac{y'_\ell}{\|y'_\ell\|}\right)
  -
  \frac{u'_{\ell+1}}{\|y'_{\ell+1}\|}\cdot\frac{y'_\ell}{\|y'_\ell\|}
  -
  \frac{u'_\ell}{\|y'_\ell\|}\cdot\frac{y'_{\ell+1}}{\|y'_{\ell+1}\|}\biggr]\notag\\
  &=
  \alpha\varepsilon\sum_{\ell=1}^N
  \biggr[(\cos{\beta_\ell}+\cos{\beta_{\ell-1}})\left(\frac{u'_\ell}{\|y'_\ell\|}\cdot\frac{y'_\ell}{\|y'_\ell\|}\right)
  -
  \frac{u'_\ell}{\|y'_\ell\|}\cdot\left(\frac{y'_{\ell+1}}{\|y'_{\ell+1}\|}+\frac{y'_{\ell-1}}{\|y'_{\ell-1}\|}\right)\biggr],\notag
\end{align}
where in the last step we rearranged the sum using the periodicity of the chain. We now notice that if all turning angles are the same, that is, if $\beta_\ell=\beta$ for all $\ell$, then we also have
\begin{equation}
  \label{eq:NNN_derivatives}
  \frac{y'_{\ell+1}}{\|y'_{\ell+1}\|}+\frac{y'_{\ell-1}}{\|y'_{\ell-1}\|}
  =
  2\cos{\beta}\frac{y'_\ell}{\|y'_\ell\|},
\end{equation}
and
\begin{equation*}
  \delta\E^{\text{b}}(y)[u]
  =
  \alpha\varepsilon\sum_{\ell=1}^N
  \biggr[2\cos{\beta}\left(\frac{u'_\ell}{\|y'_\ell\|}\cdot\frac{y'_\ell}{\|y'_\ell\|}\right)
  -
  2\cos{\beta}\left(\frac{u'_\ell}{\|y'_\ell\|}\cdot\frac{y'_\ell}{\|y'_\ell\|}\right)\biggr]
  =
  0.
\end{equation*}
This means that, among others, straight chains and uniform circular chains are critical points of this bond-angle energy.

To compute the second variation of $\E^{\text{b}}$, we start with the first expression from \eqref{eq:variation_1_Eb} and use the product rule twice. After some simplifications we obtain
\begin{align*}
  \delta^2\E^{\text{b}}(y)[u,v]
  &=
  \alpha\varepsilon\sum_{\ell=1}^N
  \biggr[\frac{u'_{\ell+1}}{\|y'_{\ell+1}\|}\cdot\left(\frac{y'_{\ell+1}}{\|y'_{\ell+1}\|}\otimes\frac{y'_\ell}{\|y'_\ell\|}+\frac{y'_\ell}{\|y'_\ell\|}\otimes\frac{y'_{\ell+1}}{\|y'_{\ell+1}\|}\right)\frac{v'_{\ell+1}}{\|y'_{\ell+1}\|}\\
  &\qquad\qquad+
  \frac{u'_\ell}{\|y'_\ell\|}\cdot\left(\frac{y'_{\ell+1}}{\|y'_{\ell+1}\|}\otimes\frac{y'_\ell}{\|y'_\ell\|}+\frac{y'_\ell}{\|y'_\ell\|}\otimes\frac{y'_{\ell+1}}{\|y'_{\ell+1}\|}\right)\frac{v'_\ell}{\|y'_\ell\|}\\
  &\qquad\qquad+
  \cos{\beta_\ell}\,\frac{u'_{\ell+1}}{\|y'_{\ell+1}\|}\cdot\left(I_2-3\,\frac{y'_{\ell+1}}{\|y'_{\ell+1}\|}\otimes\frac{y'_{\ell+1}}{\|y'_{\ell+1}\|}\right)\frac{v'_{\ell+1}}{\|y'_{\ell+1}\|}\\
  &\qquad\qquad+
  \cos{\beta_\ell}\,\frac{u'_\ell}{\|y'_\ell\|}\cdot\left(I_2-3\,\frac{y'_\ell}{\|y'_\ell\|}\otimes\frac{y'_\ell}{\|y'_\ell\|}\right)\frac{v'_\ell}{\|y'_\ell\|}.
\end{align*}
The first two and the last two terms can be combined due to the periodicity of the chain, and we obtain
\begin{align*}
  \delta^2\E^{\text{b}}(y)[u,v]
  &=
  \alpha\varepsilon\sum_{\ell=1}^N
  \biggr[\frac{u'_\ell}{\|y'_\ell\|}\cdot\left(\left(\frac{y'_{\ell+1}}{\|y'_{\ell+1}\|}+\frac{y'_{\ell-1}}{\|y'_{\ell-1}\|}\right)\otimes\frac{y'_\ell}{\|y'_\ell\|}+\frac{y'_\ell}{\|y'_\ell\|}\otimes\left(\frac{y'_{\ell+1}}{\|y'_{\ell+1}\|}+\frac{y'_{\ell-1}}{\|y'_{\ell-1}\|}\right)\right)\frac{v'_\ell}{\|y'_\ell\|}\\
  &\qquad\qquad+
  (\cos{\beta_\ell}+\cos{\beta_{\ell-1}})\,\frac{u'_\ell}{\|y'_\ell\|}\cdot\left(I_2-3\,\frac{y'_\ell}{\|y'_\ell\|}\otimes\frac{y'_\ell}{\|y'_\ell\|}\right)\frac{v'_\ell}{\|y'_\ell\|}\biggr].
\end{align*}
Finally, for a deformation $y_F$ for which all nearest neighbor interatomic distances are $F\varepsilon$ and all turning angles are $\beta$, we can again use \eqref{eq:NNN_derivatives} to get
\begin{equation}
  \label{eq:E^b''}
  \delta^2\E^{\text{b}}(y_F)[u,v]
  =
  \varepsilon\sum_{\ell=1}^N
  u'_\ell\cdot\left(\frac{2\alpha\cos{\beta}}{F^2}(I_2-P_\ell)\right)v'_\ell.
\end{equation}
Note that if the turning angles satisfy $-\pi/2\le\beta\le\pi/2$, then
\begin{equation}\label{increase}
  \delta^2\E^{\text{b}}(y_F)[u,u]
  =
  \frac{2\alpha\cos{\beta}}{F^2}\|(I_2-P_\ell)u'_\ell\|_{\ell^2_\varepsilon}^2
  \ge
  0,
\end{equation}
and if $|\beta|<\pi/2$, then
\begin{equation*}
  \delta^2\E^{\text{b}}(y_F)[u,u]
  =
  0
  \quad\text{ if and only if }\quad
  (I_2-P_\ell)u'_\ell
  =
  0
  \text{ for all }\ell,
\end{equation*}
which corresponds to pure expansion or compression of the circular chain $y_F$. Clearly, in this case the angles $\beta$ do not change and the bond-angle energy does not either.

Using the second variation \eqref{eq:E^b''} of the bond-angle energy $\E^{\text{b}}$, we can now re-state Theorems \ref{thm:stab_CB_2D}, \ref{thm:stab_a_2D}, \ref{thm:stab_CB}, and \ref{thm:stab_a} for the augmented energies $\E^{\text{CB}}(y)+\E^{\text{b}}(y)$ and $\E^{\text{a}}(y)+\E^{\text{b}}(y)$.
\begin{theorem}
  \label{thm:stability_w_bending}
  Let $\E^{\text{a,b}}(y)=\E^{\text{a}}(y)+\E^{\text{b}}(y)$ and $\E^{\text{CB,b}}(y)=\E^{\text{CB}}(y)+\E^{\text{b}}(y)$. Let $y_F\in\tilde{\Y}$ denote the 1-D configuration of atoms with nearest neighbor interatomic spacing $F\varepsilon$ so that all turning angles are $0$. Then
  \begin{equation*}
    \inf_{u\in\mathcal{U}\setminus\{0\}}
    \frac{\delta^2\E^{\text{CB,b}}(y_F)[u,u]}{\|u'\|_{\ell^2_\varepsilon}^2}
    =
    \min{\left\{\phi''(F)+4\,\phi''(2F),\,\dfrac{\phi'(F)+2\,\phi'(2F)}{F}+\dfrac{2\alpha}{F^2}\right\}},
  \end{equation*}
  and if $\phi'(2F)\ge0$ and $\phi''(2F)\le0$, then, as $\varepsilon\to0$,
  \begin{equation*}
    \inf_{u\in\mathcal{U}\setminus\{0\}}
    \frac{\delta^2\E^{\text{a,b}}(y_F)[u,u]}{\|u'\|_{\ell^2_\varepsilon}^2}
    =
    \begin{cases}
      \min{\left\{\phi''(F)+4\,\phi''(2F)+\mathcal{O}(\varepsilon^2),\,\dfrac{\phi'(F)}{F}+\dfrac{2\alpha}{F^2}\right\}} & \text{ if $N$ is even},\\
      \min{\left\{\phi''(F)+4\,\phi''(2F)+\mathcal{O}(\varepsilon^2),\,\dfrac{\phi'(F)}{F}+\dfrac{2\alpha}{F^2}+\mathcal{O}(\varepsilon)\right\}} & \text{ if $N$ is odd}.\\
    \end{cases}
  \end{equation*}

  Let $y_F\in\Y$ denote the uniform circular configuration of $N$ atoms with nearest neighbor interatomic spacing $F\varepsilon$ so that all turning angles are $\beta_\varepsilon=2\pi/N=2\pi\varepsilon$. Then, as $\varepsilon\to0$,
  \begin{equation*}
    \inf_{u\in\mathcal{U}\setminus\{0\}}
    \frac{\delta^2\E^{\text{CB,b}}(y_F)[u,u]}{\|u'\|_{\ell^2_\varepsilon}^2}
    =
    \begin{cases}
      \min{\left\{\phi''(F)+4\,\phi''(2F),\,\dfrac{\phi'(F)+2\,\phi'(2F)}{F}+\dfrac{2\alpha\cos{\beta_\varepsilon}}{F^2}\right\}} & \text{ if $N$ is even},\\
      \min{\left\{\phi''(F)+4\,\phi''(2F),\,\dfrac{\phi'(F)+2\,\phi'(2F)}{F}+\dfrac{2\alpha\cos{\beta_\varepsilon}}{F^2}+\mathcal{O}(\varepsilon)\right\}} & \text{ if $N$ is odd},
    \end{cases}
  \end{equation*}
  and if $\phi\in\mathcal{C}^3(0,\infty)$, $\phi'(2F)\ge0$, and $\phi''(2F)\le0$, then
  \begin{equation*}
    \inf_{u\in\mathcal{U}\setminus\{0\}}
    \frac{\delta^2\E^{\text{a,b}}(y_F)[u,u]}{\|u'\|_{\ell^2_\varepsilon}^2}
    =
    \begin{cases}
      \min{\left\{\phi''(F)+4\,\phi''(2F),\,\dfrac{\phi'(F)}{F}+\dfrac{2\alpha\cos{\beta_\varepsilon}}{F^2}\right\}}+\mathcal{O}(\varepsilon^2) & \text{ if $N$ is even},\\
      \min{\left\{\phi''(F)+4\,\phi''(2F),\,\dfrac{\phi'(F)}{F}+\dfrac{2\alpha\cos{\beta_\varepsilon}}{F^2}+\mathcal{O}(\varepsilon)\right\}}+\mathcal{O}(\varepsilon^2) & \text{ if $N$ is odd}.\\
    \end{cases}
  \end{equation*}
\end{theorem}
\begin{proof}
  The claims follow immediately from Theorems \ref{thm:stab_CB_2D}, \ref{thm:stab_a_2D}, \ref{thm:stab_CB}, and \ref{thm:stab_a} and the second variation \eqref{eq:E^b''} of the bond-angle energy $\E^{\text{b}}$.
\end{proof}

\begin{remark}
  Note that in Theorem \ref{thm:stability_w_bending} the expressions $\dfrac{\phi'(F)+2\,\phi'(2F)}{F}$ in the formulas for the infima of the Cauchy--Born energies and the expressions $\dfrac{\phi'(F)}{F}$ in the formulas for the infima of the atomistic energies are augmented by $\dfrac{2\alpha\cos{\beta_\varepsilon}}{F^2}>0$ if all turning angles are the same and satisfy $|\beta_\varepsilon|<\pi/2$. Also note that, for typical potentials such that $\phi'(F)<0$ if $0<F<1$, the additional bond-angle term lowers the lower bound on $F$ below which the chains can undergo the zig-zag buckling. Note, however, that for the Lennard-Jones or Morse potentials, the term $\dfrac{\phi'(F)}{F}$ is negative and dominates the term $\dfrac{2\alpha\cos{\beta_\varepsilon}}{F^2}$ as $F\to0$, so no matter how large $\alpha$ is, there is always a positive lower bound for the region of stability of the deformation $y_F$.
\end{remark}

\begin{remark}
  The bond-angle energy \eqref{eq:E^b} given by $\E^{\text{b}}(y)=\varepsilon\sum_{\ell=1}^N\alpha(1-\cos{\beta_\ell})$ is meant to model in one-dimensional chains the resistance to transverse displacement in two-dimensional graphene and carbon nanotubes. There is an energy cost to transverse displacement in tri-bonded graphene and carbon nanotubes when modeled by popular potentials such as \cite{tersoff88} which penalize the bond-angle deviation from $2\pi/3$.
\end{remark}

\section{Modeling errors}
\label{sec:truncation}
In this section, we will consider applying an external periodic load $f\in\U$ to a deformation $y_F\in\Y$ and study the error of the Cauchy--Born model. More specifically, let $u^{\text{a}}\in\U$ and $u^{\text{CB}}\in\U$ solve the linearized equations
\begin{alignat*}
  \delta\E^{\text{a}}(y_F)[v]+\delta^2\E^{\text{a}}(y_F)[u^{\text{a}},v]
  &=
  \langle f,v\rangle
  &\qquad\text{ for all }v\in\U,\\
  \delta\E^{\text{CB}}(y_F)[v]+\delta^2\E^{\text{CB}}(y_F)[u^{\text{CB}},v]
  &=
  \langle f,v\rangle
  &\qquad\text{ for all }v\in\U.
\end{alignat*}
We define the modeling error of the Cauchy--Born approximation, $\tau$, via the duality relationship
\begin{equation}
  \label{eq:trunc_def}
  \begin{split}
  \langle\tau,v\rangle
  :&=
  \delta\E^{\text{CB}}(y_F)[v]+\delta^2\E^{\text{CB}}(y_F)[u^{\text{a}},v]-\langle f,v\rangle\\
  &=
  \left\{ \delta\E^{\text{CB}}(y_F)[v]-\delta\E^{\text{a}}(y_F)[v]\right\}
  +
  \left\{ \delta^2\E^{\text{CB}}(y_F)[u^{\text{a}},v]-\delta^2\E^{\text{a}}(y_F)[u^{\text{a}},v]\right\}
  \quad\text{ for all }v\in\U.
  \end{split}
\end{equation}
Since the solution $u^{\text{CB}}$ will play no role in the analysis of the modeling error, to simplify the notation in the rest of this section, we will sometimes suppress the superscript and simply write $u$ instead of $u^{\text{a}}$. However, in the statements of the theorems, the proper notation will be used.

We compute that we have for all $v\in\U$
\begin{equation}
  \label{eq:ghost}
  \delta\E^{\text{CB}}(y_F)[v]-\delta\E^{\text{a}}(y_F)[v]
  =
  \varepsilon\sum_{\ell=1}^N\left[\left(\frac{\phi'(2F_1)}{2F_1}-\frac{\phi'(2F_2)}{2F_2}\right)\left(4y'_{F,\ell}+\varepsilon^2y'''_{F,\ell+1}\right)-\varepsilon^2\frac{\phi'(2F_1)}{2F_1}y'''_{F,\ell+1}\right]\cdot v'_\ell
\end{equation}
and thus observe that
\begin{equation*}
  \delta\E^{\text{CB}}(y_F)-\delta\E^{\text{a}}(y_F)\equiv0
  \qquad\text{if $y_F$ is a linear chain since $F_1=F_2$ and $y'''_{F}\equiv0$,}
\end{equation*}
but that
\begin{equation*}
  \delta\E^{\text{CB}}(y_F)-\delta\E^{\text{a}}(y_F)\not\equiv0
  \qquad\text{if $y_F$ is a circular chain.}
\end{equation*}

\subsection{Modeling errors for linear chains}
Using the second variations of $\E^{\text{a}}$ and $\E^{\text{CB}}$ given in \eqref{eq:E^a''_second} and \eqref{eq:E^CB''}, respectively, we now have since $\delta\E^{\text{CB}}(y_F)-\delta\E^{\text{a}}(y_F)\equiv 0$ for linear chains that
\begin{align*}
  \langle\tau,v\rangle
  &=
  \delta^2\E^{\text{CB}}(y_F)[u,v]
  -
  \delta^2\E^{\text{a}}(y_F)[u,v]\\
  &=
  \varepsilon\sum_{\ell=1}^N
  \biggr[
  u'_{\ell}\cdot\left(\bigr(\phi''(F_1)+4\,\phi''(2F_1)\bigr)P_\ell+\frac{\phi'(F_1)+2\,\phi'(2F_1)}{F_1}\left(I_2-P_\ell\right)\right)v'_{\ell}\\
  &\qquad\qquad-
  u'_{\ell}\cdot\left(\phi''(F_1)P_\ell+\frac{\phi'(F_1)}{F_1}(I_2-P_\ell)\right)v'_{\ell}\\
  &\qquad\qquad-
  (u'_{\ell+1}+u'_{\ell})\cdot\left(\phi''(2F_2)\tilde{P}_\ell+\frac{\phi'(2F_2)}{2F_2}(I_2-\tilde{P}_\ell)\right)(v'_{\ell+1}+v'_\ell)\biggr].
\end{align*}
If we define the matrices
\begin{equation}
  \label{eq:A_l}
  A_\ell
  :=
  \phi''(2F_1)P_\ell+\frac{\phi'(2F_1)}{2F_1}(I_2-P_\ell)
  \quad\text{ and }\quad
  \tilde{A}_\ell
  :=
  \phi''(2F_2)\tilde{P}_\ell+\frac{\phi'(2F_2)}{2F_2}(I_2-\tilde{P}_\ell),
\end{equation}
cancel the terms on the second line and use identity \eqref{eq:parallelogram}, after rearranging the sum we get
\begin{align*}
  \langle\tau,v\rangle
  &=
  \varepsilon\sum_{\ell=1}^N
  \biggr[
  4\,u'_{\ell}\cdot A_\ell v'_{\ell}
  -
  2\,u'_{\ell+1}\tilde{A}_\ell v'_{\ell+1}
  -
  2\,u'_\ell\tilde{A}_\ell v'_\ell
  +
  \varepsilon^2u''_{\ell+1}\tilde{A}_\ell v''_{\ell+1}\biggr]\\
  &=
  \varepsilon\sum_{\ell=1}^N
  \biggr[
  \varepsilon^2u''_{\ell+1}\cdot\tilde{A}_\ell v''_{\ell+1}
  -
  2\,u'_\ell\cdot\left(\tilde{A}_\ell+\tilde{A}_{\ell-1}-2A_\ell\right)v'_\ell\biggr].
\end{align*}
We can now sum by parts and simplify to get
\begin{align}
  \label{eq:trunc}
  \langle\tau,v\rangle
  &=
  -\varepsilon\sum_{\ell=1}^N
  \biggr[\varepsilon^2\left(\tilde{A}_\ell u'''_{\ell+1}+\frac{\tilde{A}_\ell-\tilde{A}_{\ell-1}}{\varepsilon}u''_\ell\right)+
  2\left(\tilde{A}_\ell+\tilde{A}_{\ell-1}-2A_\ell\right)u'_\ell\biggr]\cdot v'_\ell.
\end{align}

Let us now consider a 1-D atomic configuration $y_F\in\tilde{\Y}$ with interatomic spacing $F\varepsilon$. In this case we have $F_1=F_2=F$, and also all of the projection operators $P_\ell$ and $\tilde{P}_\ell$ are the same for all $\ell$, so we can write $P_\ell=\tilde{P}_\ell=P$ and $A_\ell=\tilde{A}_\ell=A$. The expression \eqref{eq:trunc} for the modeling error then simplifies to
\begin{equation*}
  \langle\tau,v\rangle
  =
  -\varepsilon\sum_{\ell=1}^N
  \varepsilon^2(Au'''_{\ell+1})\cdot v'_\ell,
\end{equation*}
and we immediately have the following theorems, one for the 1-D constrained chain and one for the 1-D unconstrained chain.
\begin{theorem}
  \label{thm:truncation_1D1D}
  Let $y_F\in\tilde{\Y}$ denote the 1-D configuration of atoms with nearest neighbor interatomic spacing $F\varepsilon$, let $f\in\tilde{\U}$, and let $u^{\text{a}}\in\tilde{\U}$ satisfy $\delta\E^{\text{a}}(y_F)[v]+\delta^2\E^{\text{a}}(y_F)[u^{\text{a}},v]=\langle f,v\rangle$ for all $v\in\tilde{\U}$. The modeling error of the Cauchy--Born approximation, $\tau$, then satisfies the inequality
  \begin{equation}
    \label{eq:trunc_1D1D}
    \|\tau\|_*
    \le
    \varepsilon^2|\phi''(2F)|\,\|(u^{\text{a}})'''\|_{\ell_\varepsilon^2}.
  \end{equation}
\end{theorem}
\begin{theorem}
  \label{thm:truncation_1D2D}
  Let $y_F\in\tilde{\Y}$ denote the 1-D configuration of atoms with nearest neighbor interatomic spacing $F\varepsilon$, let $f\in\U$, and let $u^{\text{a}}\in\U$ satisfy $\delta\E^{\text{a}}(y_F)[v]+\delta^2\E^{\text{a}}(y_F)[u^{\text{a}},v]=\langle f,v\rangle$ for all $v\in\U$. The modeling error of the Cauchy--Born approximation, $\tau$, then satisfies the inequality
  \begin{equation}
    \label{eq:trunc_1D2D}
    \|\tau\|_*
    \le
    \varepsilon^2\max{\left\{|\phi''(2F)|,\biggr|\frac{\phi'(2F)}{2F}\biggr|\right\}}\|(u^{\text{a}})'''\|_{\ell_\varepsilon^2}.
  \end{equation}
\end{theorem}

\subsection{Modeling error for circular chains}
The situation is different in the case of a uniform circular chain $y_F\in\Y$ with interatomic spacing $F\varepsilon$. In this case, $\delta\E^{\text{CB}}(y_F)-\delta\E^{\text{a}}(y_F)\not\equiv0$, but recalling the constant $C_\phi$ from \eqref{eq:C_phi}, we obtain from \eqref{eq:ghost}
\begin{align}
  \label{eq:ghost_estimate}
    \left|\delta\E^{\text{CB}}(y_F)[v]-\delta\E^{\text{a}}(y_F)[v]\right|
    &=
    \left|\varepsilon\sum_{\ell=1}^N\left[\left(\frac{\phi'(2F_1)}{2F_1}-\frac{\phi'(2F_2)}{2F_2}\right)\left(4y'_{F,\ell}+\varepsilon^2y'''_{F,\ell+1}\right)-\varepsilon^2\frac{\phi'(2F_1)}{2F_1}y'''_{F,\ell+1}\right]\cdot v'_\ell\right|\notag\\
    &\le
    \left[\varepsilon^2C_\phi\left(4\|y'_F\|_{\ell^2_\varepsilon}+\varepsilon^2\|y'''_F\|_{\ell^2_\varepsilon}\right)+\varepsilon^2\left|\frac{\phi'(2F)}{2F}\right|\|y'''_F\|_{\ell^2_\varepsilon}\right]\|v'\|_{\ell^2_\varepsilon}\notag\\
    &\le
    \varepsilon^2\biggr[4\,C_\phi(1+\pi^2\varepsilon^2)F+2\pi^2|\phi'(2F)|\biggr]\|v'\|_{\ell^2_\varepsilon}\\
    &=
    \varepsilon^2C_\kappa\|v'\|_{\ell^2_\varepsilon},\notag
\end{align}
where
\begin{equation}
  \label{eq:C_kappa}
    C_\kappa
    =
    4\,C_\phi(1+\pi^2\varepsilon^2)F+2\pi^2|\phi'(2F)|
\end{equation}
since $\|y'_F\|_{\ell^2_\varepsilon}=F$ and $\|y'''_F\|_{\ell^2_\varepsilon}=4F\varepsilon^{-2}\sin^2(\pi\varepsilon)\le4F\pi^2$.

For the modeling error term $\delta^2\E^{\text{CB}}(y_F)[u,v]-\delta^2\E^{\text{a}}(y_F)[u,v]$, there is no cancellation in \eqref{eq:trunc} as in the case of the straight chain, but the expressions $\tilde{A}_\ell-\tilde{A}_{\ell-1}$ and $\tilde{A}_\ell+\tilde{A}_{\ell-1}-2A_\ell$ are of order $\varepsilon$, and $\varepsilon^2$, respectively. Therefore, the modeling error of the Cauchy--Born approximation is again
of order $\varepsilon^2$.
\begin{theorem}
  \label{thm:truncation_2D}
  Let $y_F\in\Y$ denote the uniform circular configuration of $N$ atoms with nearest neighbor interatomic spacing $F\varepsilon$, let $f\in\U$, and let $u^{\text{a}}\in\U$ satisfy $\delta\E^{\text{a}}(y_F)[v]+\delta^2\E^{\text{a}}(y_F)[u^{\text{a}},v]=\langle f,v\rangle$ for all $v\in\U$. If $\phi\in\C^3(0,\infty)$, then the modeling error of the Cauchy--Born approximation, $\tau$, satisfies the inequality
  \begin{equation}
    \label{eq:trunc_2D}
    \|\tau\|_*
    \le
    C_\kappa \varepsilon^2
    +
    (C_1\varepsilon^2+C_2\varepsilon^4)\left(\|(u^{\text{a}})'''\|_{\ell_\varepsilon^2}+\|(u^{\text{a}})''\|_{\ell_\varepsilon^2}
    +\|(u^{\text{a}})'\|_{\ell_\varepsilon^2}\right),
  \end{equation}
  where
  \begin{equation}
    \label{eq:C1C2}
    C_1
    =
    \max\left\{|\phi''(2F)|,\ \biggr|\frac{\phi'(2F)}{2F}\biggr|,\ 4\pi^2\left|\phi''(2F)-\frac{\phi'(2F)}{2F}\right|+12\,C_\phi\right\},
    \quad
    C_2
    =
    4\pi\,C_\phi,
  \end{equation}
  $C_\phi$ is the Lipschitz constant defined in \eqref{eq:C_phi}, and $C_\kappa$ is defined in \eqref{eq:C_kappa}.
\end{theorem}
\begin{proof}
  Using the definition \eqref{eq:A_l} of $A_\ell$ and $\tilde{A}_\ell$, we have, after some manipulations,
  \begin{equation*}
    \tilde{A}_\ell-\tilde{A}_{\ell-1}
    =
    \left(\phi''(2F_2)-\frac{\phi'(2F_2)}{2F_2}\right)\left(\tilde{P}_\ell-\tilde{P}_{\ell-1}\right),
  \end{equation*}
  and
  \begin{align*}
    \tilde{A}_\ell+\tilde{A}_{\ell-1}-2A_\ell
    &=
    \left(\phi''(2F_1)-\frac{\phi'(2F_1)}{2F_1}\right)\left(\tilde{P}_\ell+\tilde{P}_{\ell-1}-2P_\ell\right)
    +
    2\left(\frac{\phi'(2F_2)}{2F_2}-\frac{\phi'(2F_1)}{2F_1}\right)I_2\\
    &\quad\quad+
    \left[\left(\phi''(2F_2)-\phi''(2F_1)\right)-\left(\frac{\phi'(2F_2)}{2F_2}-\frac{\phi'(2F_1)}{2F_1}\right)\right]\left(\tilde{P}_\ell+\tilde{P}_{\ell-1}\right).
  \end{align*}
  It can be easily verified that for any vector $w\in\Real^2$ we have (cf.~\eqref{eq:projections_2})
  \begin{gather*}
    \|(\tilde{P}_\ell-\tilde{P}_{\ell-1})w\|=\|w\|\sin\frac{2\pi}{N}=\|w\|\sin{2\pi\varepsilon}\le2\pi\varepsilon\|w\|,\\
    \|(\tilde{P}_\ell+\tilde{P}_{\ell-1}-2P_\ell)w\|=2\|w\|\sin^2\frac{\pi}{N}=2\|w\|\sin^2{(\pi\varepsilon)}\le2\pi^2\varepsilon^2\|w\|.
  \end{gather*}
  Recall that $F_1=F$ and $F_2=F\cos{(\pi\varepsilon)}$, and from \eqref{eq:C_phi} that there exists a constant $C_\phi>0$ independent of $\varepsilon$ such that
  \begin{equation*}
    \max\left\{\left|\phi''(2F_2)-\phi''(2F_1)\right|,\left|\frac{\phi'(2F_2)}{2F_2}-\frac{\phi'(2F_1)}{2F_1}\right|\right\}
    \le
    C_\phi\,\varepsilon^2.
  \end{equation*}
  Therefore, from \eqref{eq:trunc} and \eqref{eq:ghost_estimate} we obtain
  \begin{align*}
    \|\tau\|_*
    &\le
    \varepsilon^2\max{\left\{|\phi''(2F_2)|,\biggr|\frac{\phi'(2F_2)}{2F_2}\biggr|\right\}}\|(u^{\text{a}})'''\|_{\ell_\varepsilon^2}
    +
    2\pi\varepsilon^2\left|\phi''(2F_2)-\frac{\phi'(2F_2)}{2F_2}\right|\|(u^{\text{a}})''\|_{\ell^2_\varepsilon}\\
    &\quad\quad+
    2\left(2\pi^2\varepsilon^2\left|\phi''(2F_1)-\frac{\phi'(2F_1)}{2F_1}\right|+6\,C_\phi\,\varepsilon^2\right)\|(u^{\text{a}})'\|_{\ell_\varepsilon^2}
    +
    C_\kappa\varepsilon^2\\
    &\le
    \varepsilon^2\left(\max{\left\{|\phi''(2F_1)|,\biggr|\frac{\phi'(2F_1)}{2F_1}\biggr|\right\}}+C_\phi\,\varepsilon^2\right)\|(u^{\text{a}})'''\|_{\ell_\varepsilon^2}\\
    &\qquad+
    \varepsilon^2\left(2\pi\left|\phi''(2F_1)-\frac{\phi'(2F_1)}{2F_1}\right|+4\pi\,C_\phi\,\varepsilon^2\right)\|(u^{\text{a}})''\|_{\ell^2_\varepsilon}\\
    &\quad\quad+
    \varepsilon^2\left(4\pi^2\left|\phi''(2F_1)-\frac{\phi'(2F_1)}{2F_1}\right|+12\,C_\phi\right)\|(u^{\text{a}})'\|_{\ell_\varepsilon^2}
    +
    C_\kappa\varepsilon^2
  \end{align*}
  and the statement of the theorem follows.
\end{proof}
\begin{remark}
 The above second-order modeling error estimates $\mathcal{O}(\varepsilon^{2})$ for the Cauchy--Born approximation were the result of the symmetric treatment of the interactions between next nearest neighbors. A more explicit treatment of the second-order modeling error for the Cauchy--Born approximation for linear chains was given in~\cite{Dobson:2008b}.

 We note, however, that the second-order modeling error estimates $\mathcal{O}(\varepsilon^{2})$ for the Cauchy--Born approximation require that $\|(u^{\text{a}})'''\|_{\ell_\varepsilon^2}$ be bounded uniformly in $\varepsilon$, which is not the case for the approximation of atomistic configurations with defects. This lack of accuracy of the Cauchy--Born approximation for problems with defects is the motivation for the development of atomistic-to-continuum methods such as the quasi-nonlocal method which attain $\mathcal{O}(\varepsilon^{3/2})$ accuracy for problems with defects~\cite{dobs-qcf2,brian10,ortner:qnl1d}.
\end{remark}

\section{Error analysis for the Cauchy--Born approximation}
\label{sec:error}
From the definition \eqref{eq:trunc_def} of the modeling error of the Cauchy--Born approximation, we have
\begin{equation}
  \label{eq:trunc_def2}
  \begin{split}
  \langle\tau,v\rangle
  &=
  \delta\E^{\text{CB}}(y_F)[v]+\delta^2\E^{\text{CB}}(y_F)[u^{\text{a}},v]-\langle f,v\rangle\\
  &=
  \delta\E^{\text{CB}}(y_F)[v]+\delta^2\E^{\text{CB}}(y_F)[u^{\text{a}},v]
  -
  \delta\E^{\text{CB}}(y_F)[v]-\delta^2\E^{\text{CB}}(y_F)[u^{\text{CB}},v]\\
  &=
  \delta^2\E^{\text{CB}}(y_F)[u^{\text{a}},v]
  -
  \delta^2\E^{\text{CB}}(y_F)[u^{\text{CB}},v]\\
  &=
  \delta^2\E^{\text{CB}}(y_F)[u^{\text{a}}-u^{\text{CB}},v]
  \qquad\text{ for all }v\in\U.
  \end{split}
\end{equation}
Setting $v=u^{\text{a}}-u^{\text{CB}}$ in \eqref{eq:trunc_def2} above, we obtain that
\begin{equation}
  \label{eq:trunc_stab}
  \delta^2\mathcal{E}^{\text{CB}}(y_F)[u^{\text{a}}-u^{\text{CB}},u^{\text{a}}-u^{\text{CB}}]
  =
  \langle\tau,u^{\text{a}}-u^{\text{CB}}\rangle.
\end{equation}

Let us also define the Cauchy--Born stability constants (see
Theorems \ref{thm:stab_CB_1D}, \ref{thm:stab_CB_2D}, and
\ref{thm:stab_CB})
\begin{equation*}
  \gamma_1
  =
  \phi''(F)+4\,\phi''(2F)
  \quad\text{ and }\quad
  \gamma_2
  =
  \min{\left\{\phi''(F)+4\,\phi''(2F),\,\dfrac{\phi'(F)+2\,\phi'(2F)}{F}\right\}}.
\end{equation*}
Using the stability results and the expressions for the modeling errors, we now have the following theorems for the errors of the Cauchy--Born approximations.

\begin{theorem}
  \label{thm:error_1D1D}
  Let $y_F\in\tilde{\Y}$ denote the 1-D configuration of atoms with nearest neighbor interatomic spacing $F\varepsilon$. Given $f\in\tilde{\U}$, let $u^{\text{a}},u^{\text{CB}}\in\tilde{\U}$ satisfy $\delta\E^{\text{a}}(y_F)[v]+\delta^2\E^{\text{a}}(y_F)[u^{\text{a}},v]=\langle f,v\rangle$ and $\delta\E^{\text{CB}}(y_F)[v]+\delta^2\E^{\text{CB}}(y_F)[u^{\text{CB}},v]=\langle f,v\rangle$, respectively, for all $v\in\tilde{\U}$. If $\gamma_1>0$, then
  \begin{equation*}
    \|(u^{\text{a}}-u^{\text{CB}})'\|_{\ell^2_\varepsilon}
    \le \frac{\varepsilon^2 |\phi''(2F)|}{\gamma_1}\|(u^{\text{a}})'''\|_{\ell_\varepsilon^2}.
  \end{equation*}
\end{theorem}
\begin{proof}
  Combining \eqref{eq:trunc_stab} with \eqref{eq:stab_CB_1D1D_estimate} and \eqref{eq:trunc_1D1D}, the theorem follows from the inequality
  \begin{align*}
    \gamma_1\|(u^{\text{a}}-u^{\text{CB}})'\|^2_{\ell^2_\varepsilon}
    &\le
    \delta^2\E^{\text{CB}}(y_F)[u^{\text{a}}-u^{\text{CB}},u^{\text{a}}-u^{\text{CB}}]\\
    &=
    \langle\tau,u^{\text{a}}-u^{\text{CB}}\rangle\\
    &\le
    \|\tau\|_*\|(u^{\text{a}}-u^{\text{CB}})'\|_{\ell^2_\varepsilon}\\
    &\le
    \varepsilon^2|\phi''(2F)|\,\|(u^{\text{a}})'''\|_{\ell_\varepsilon^2}
    \left\|(u^{\text{a}}-u^{\text{CB}})'\right\|_{\ell^2_\varepsilon}.
  \end{align*}
\end{proof}

\begin{theorem}
  \label{thm:error_1D2D}
  Let $y_F\in\tilde{\Y}$ denote the 1-D configuration of atoms with nearest neighbor interatomic spacing $F\varepsilon$. Given $f\in\U$, let $u^{\text{a}},u^{\text{CB}}\in\U$ satisfy $\delta\E^{\text{a}}(y_F)[v]+\delta^2\E^{\text{a}}(y_F)[u^{\text{a}},v]=\langle f,v\rangle$ and $\delta\E^{\text{CB}}(y_F)[v]+\delta^2\E^{\text{CB}}(y_F)[u^{\text{CB}},v]=\langle f,v\rangle$, respectively, for all $v\in\U$. If $\gamma_2>0$, then
  \begin{equation*}
    \|(u^{\text{a}}-u^{\text{CB}})'\|_{\ell^2_\varepsilon}
    \le
    \frac{\varepsilon^2\max{\left\{|\phi''(2F)|,\biggr|\dfrac{\phi'(2F)}{2F}\biggr|\right\}}}{\gamma_2}
    \left\|(u^{\text{a}})'''\right\|_{\ell_\varepsilon^2}.
  \end{equation*}
\end{theorem}
\begin{proof}
  Combining \eqref{eq:trunc_stab} with \eqref{eq:stab_CB_1D2D_estimate} and \eqref{eq:trunc_1D2D}, the theorem follows from the inequality
  \begin{align*}
    \gamma_2\|(u^{\text{a}}-u^{\text{CB}})'\|^2_{\ell^2_\varepsilon}
    &\le
    \delta^2\E^{\text{CB}}(y_F)[u^{\text{a}}-u^{\text{CB}},u^{\text{a}}-u^{\text{CB}}]\\
    &=
    \langle\tau,u^{\text{a}}-u^{\text{CB}}\rangle\\
    &\le
    \|\tau\|_*\|(u^{\text{a}}-u^{\text{CB}})'\|_{\ell^2_\varepsilon}\\
    &\le
    \varepsilon^2\max{\left\{|\phi''(2F)|,\biggr|\frac{\phi'(2F)}{2F}\biggr|\right\}}\|(u^{\text{a}})'''\|_{\ell_\varepsilon^2}\|(u^{\text{a}}-u^{\text{CB}})'\|_{\ell^2_\varepsilon}.
  \end{align*}
\end{proof}

\begin{theorem}
  \label{thm:error_2D}
  Let $y_F\in\Y$ denote the uniform circular configuration of $N$ atoms with nearest neighbor interatomic spacing $F\varepsilon$. Given $f\in\U$, let $u^{\text{a}},u^{\text{CB}}\in\U$ satisfy $\delta\E^{\text{a}}(y_F)[v]+\delta^2\E^{\text{a}}(y_F)[u^{\text{a}},v]=\langle f,v\rangle$ and $\delta\E^{\text{CB}}(y_F)[v]+\delta^2\E^{\text{CB}}(y_F)[u^{\text{CB}},v]=\langle f,v\rangle$, respectively, for all $v\in\U$. If $\phi\in\mathcal{C}^3(0,\infty)$ and $\gamma_2>0$, then
  \begin{equation*}
    \|(u^{\text{a}}-u^{\text{CB}})'\|_{\ell^2_\varepsilon}
    \le
    \gamma_2^{-1}\left[C_\kappa\varepsilon^2+(C_1\varepsilon^2+C_2\varepsilon^4)\left(\|(u^{\text{a}})'''\|_{\ell_\varepsilon^2}+\|(u^{\text{a}})''\|_{\ell_\varepsilon^2}+\|(u^{\text{a}})'\|_{\ell_\varepsilon^2}\right)\right],
  \end{equation*}
  where $C_1$ and $C_2$ are defined in \eqref{eq:C1C2} and $C_\kappa$ in \eqref{eq:C_kappa}.
\end{theorem}
\begin{proof}
  Combining \eqref{eq:trunc_stab} with \eqref{eq:stab_CB_2D_estimate} and \eqref{eq:trunc_2D}, we obtain
  \begin{align*}
    \gamma_2\|(u^{\text{a}}-u^{\text{CB}})'\|^2_{\ell^2_\varepsilon}
    &\le
    \delta^2\E^{\text{CB}}(y_F)[u^{\text{a}}-u^{\text{CB}},u^{\text{a}}-u^{\text{CB}}]\\
    &=
    \langle\tau,u^{\text{a}}-u^{\text{CB}}\rangle\\
    &\le
    \|\tau\|_*\|(u^{\text{a}}-u^{\text{CB}})'\|_{\ell^2_\varepsilon}\\
    &\le
    \left[C_\kappa\varepsilon^2+(C_1\varepsilon^2+C_2\varepsilon^4)\left(\|(u^{\text{a}})'''\|_{\ell_\varepsilon^2}+\|(u^{\text{a}})''\|_{\ell_\varepsilon^2}+\|(u^{\text{a}})'\|_{\ell_\varepsilon^2}\right)\right]\|(u^{\text{a}}-u^{\text{CB}})'\|_{\ell^2_\varepsilon},
  \end{align*}
  and the theorem follows.
\end{proof}

\begin{remark}
  We note that the error $\|(u^{\text{a}}-u^{\text{CB}})'\|_{\ell^2_\varepsilon}$ will generally not be small for compressive $F<1$ such that the atomistic model is not stable ($\phi'(F)<0$), but the Cauchy--Born model is stable ($\phi'(F)+2\phi'(2F)>0$), because $\|(u^{\text{a}})'''\|_{\ell_\varepsilon^2}$ can be expected to be large.
\end{remark}

\section{Error analysis of the models including the bond-angle energy}
\label{sec:error_bending}
The error analysis for the case when the bond-angle energy $\E^{\text{b}}$ as defined in \eqref{eq:E^b} is included in the models is straightforward. The stability results are given in Theorem \ref{thm:stability_w_bending}. Since we only need an estimate from below on the Cauchy--Born stability constant, we can use again $\gamma_2$ as defined in Section \ref{sec:error},
\begin{equation*}
  \gamma_2
  =
  \min{\left\{\phi''(F)+4\,\phi''(2F),\,\dfrac{\phi'(F)+2\,\phi'(2F)}{F}\right\}}.
\end{equation*}
There is no contribution of the bond-angle energy $\E^{\text{b}}$ to the modeling error, $\tau$, since  $\tau$ is given in \eqref{eq:trunc_def} by
\begin{align*}
  \langle\tau,v\rangle
  &=
  \left\{ \delta\E^{\text{CB,b}}(y_F)[v]-\delta\E^{\text{a,b}}(y_F)[v]\right\}
  +
  \left\{\delta^2\E^{\text{CB,b}}(y_F)[u^{\text{a,b}},v]-\delta^2\E^{\text{a,b}}(y_F)[u^{\text{a,b}},v]\right\}\\
  &=
  \left\{ \delta\E^{\text{CB}}(y_F)[v]-\delta\E^{\text{a}}(y_F)[v]\right\}
  +
  \left\{ \delta^2\E^{\text{CB}}(y_F)[u^{\text{a,b}},v]-\delta^2\E^{\text{a}}(y_F)[u^{\text{a,b}},v]\right\},
\end{align*}
where $\delta\E^{\text{a,b}}(y_F)[v]+\delta^2\E^{\text{a,b}}(y_F)[u^{\text{a,b}},v]=\langle f,v\rangle$ and $\delta\E^{\text{CB,b}}(y_F)[v]+\delta^2\E^{\text{CB,b}}(y_F)[u^{\text{CB,b}},v]=\langle f,v\rangle$ for all $v\in\U$. Therefore, the estimates \eqref{eq:trunc_1D1D}, \eqref{eq:trunc_1D2D}, and \eqref{eq:trunc_2D} for the modeling error still hold with $u^{\text{a}}$ replaced by $u^{\text{a,b}}$.

We therefore have the following theorems, which are analogous to Theorems \ref{thm:error_1D2D} and \ref{thm:error_2D}.
\begin{theorem}
  Let $y_F\in\tilde{\Y}$ denote the 1-D configuration of atoms with nearest neighbor interatomic spacing $F\varepsilon$. Given $f\in\U$, let $u^{\text{a,b}},u^{\text{CB,b}}\in\U$ satisfy $\delta\E^{\text{a,b}}(y_F)[v]+\delta^2\E^{\text{a,b}}(y_F)[u^{\text{a,b}},v]=\langle f,v\rangle$ and $\delta\E^{\text{CB,b}}(y_F)[v]+\delta^2\E^{\text{CB,b}}(y_F)[u^{\text{CB,b}},v]=\langle f,v\rangle$, respectively, for all $v\in\U$. If $\gamma_2>0$, then
  \begin{equation*}
    \|(u^{\text{a,b}}-u^{\text{CB,b}})'\|_{\ell^2_\varepsilon}
    \le
   \frac{\varepsilon^2\max{\left\{|\phi''(2F)|,\biggr|\dfrac{\phi'(2F)}{2F}\biggr|\right\}}}{\gamma_2}\|(u^{\text{a,b}})'''\|_{\ell_\varepsilon^2}.
  \end{equation*}
\end{theorem}
\begin{theorem}
Let $y_F\in\Y$ denote the uniform circular configuration of $N$
atoms with nearest neighbor interatomic spacing $F\varepsilon.$
Given $f\in\U$, let
$u^{\text{a,b}},u^{\text{CB,b}}\in\U$ satisfy
$\delta\E^{\text{a,b}}(y_F)[v]+\delta^2\E^{\text{a,b}}(y_F)[u^{\text{a,b}},v]=\langle
f,v\rangle$ and
$\delta\E^{\text{CB,b}}(y_F)[v]+\delta^2\E^{\text{CB,b}}(y_F)[u^{\text{CB,b}},v]=\langle
f,v\rangle$, respectively, for all $v\in\U$. If
$\phi\in\mathcal{C}^3(0,\infty)$ and $\gamma_2>0$, then
  \begin{equation*}
    \|(u^{\text{a,b}}-u^{\text{CB,b}})'\|_{\ell^2_\varepsilon}
    \le
    \gamma_2^{-1}\left[C_\kappa\varepsilon^2+(C_1\varepsilon^2+C_2\varepsilon^4)
    \left(\|(u^{\text{a,b}})'''\|_{\ell_\varepsilon^2}+\|(u^{\text{a,b}})''\|_{\ell_\varepsilon^2}
    +
    \|(u^{\text{a,b}})'\|_{\ell_\varepsilon^2}\right)\right],
  \end{equation*}
  where $C_1$ and $C_2$ are defined in \eqref{eq:C1C2} and $C_\kappa$ in \eqref{eq:C_kappa}.
\end{theorem}

\section{The quasi-nonlocal approximation}
\label{sec:qnl}
In this section, we consider a combination of the atomistic approach in one region and the Cauchy--Born approximation in its complement. Specifically, for an integer $1<K<N$, we define the set $\mathfrak{A}$ of \emph{nonlocal} atoms (``nonlocal'' or ``atomistic'' region) and the set $\mathfrak{C}$ of \emph{local} atoms (``local'' or ``continuum'' region) as
\begin{equation*}
  \mathfrak{A}
  =
  \{1,\dots,K\}
  \quad\text{ and }\quad
  \mathfrak{C}
  =
  \{K+1,\dots,N\},
\end{equation*}
extended by periodicity.

For $y\in\Y$, we define the quasi-nonlocal energy, $\E^{\text{QNL}}(y)$, as a combination of the atomistic energy in the nonlocal region and the Cauchy--Born approximation in the local region:
\begin{align*}
  \E^{\text{QNL}}(y)
  &=
  \varepsilon\sum_{\ell\in\mathfrak{A}}
  \biggr[\phi(\|y'_\ell\|)+\phi(\|y'_{\ell+1}+y'_\ell\|)\biggr]
  +
  \varepsilon\sum_{\ell\in\mathfrak{C}}
  \biggr[\phi(\|y'_\ell\|)+\frac12\left(\phi(2\|y'_\ell\|)+\phi(2\|y'_{\ell+1}\|)\right)\biggr]\\
  &=
  \varepsilon\sum_{\ell=1}^N\phi(\|y'_\ell\|)
  +
  \varepsilon\sum_{\ell=1}^K\phi(\|y'_{\ell+1}+y'_\ell\|)
  +
  \varepsilon\sum_{\ell=K+1}^N\frac12\left(\phi(2\|y'_\ell\|)+\phi(2\|y'_{\ell+1}\|)\right).
\end{align*}
We can rearrange the sums to resemble the atomistic and Cauchy--Born energies in the following way:
\begin{equation*}
  \E^{\text{QNL}}(y)
  =
  \varepsilon\sum_{\ell=1}^N\phi(\|y'_\ell\|)
  +
  \varepsilon\sum_{\ell=1}^K\phi(\|y'_{\ell+1}+y'_\ell\|)
  +
  \varepsilon\sum_{\ell=K+2}^N\phi(2\|y'_\ell\|)
  +
  \frac{\varepsilon}{2}\,\phi(2\|y'_1\|)
  +
  \frac{\varepsilon}{2}\,\phi(2\|y'_{K+1}\|).
\end{equation*}

When studying the quasi-nonlocal method and its properties, the following seminorms of displacements $u\in\U$ will be used:
\begin{equation*}
  \|u'\|_{\ell^2_\varepsilon(\mathfrak{A})}
  =
  \left(\varepsilon\sum_{\ell=2}^K\|u'_\ell\|^2\right)^{1/2},
  \quad
  \|u'\|_{\ell^2_\varepsilon(\mathfrak{C})}
  =
  \left(\varepsilon\sum_{\ell=K+2}^N\|u'_\ell\|^2\right)^{1/2},
  \quad
  \|u'\|_{\ell^2_\varepsilon(\mathfrak{I})}
  =
  \biggr(\varepsilon\|u'_1\|^2+\varepsilon\|u'_{K+1}\|^2\biggr)^{1/2}.
\end{equation*}
Note that since the primes denote backward differences, $\|u'\|_{\ell^2_\varepsilon(\mathfrak{A})}$ is a seminorm over the bonds between nonlocal atoms, $\|u'\|_{\ell^2_\varepsilon(\mathfrak{C})}$ is a seminorm over the bonds between local atoms, and $\|u'\|_{\ell^2_\varepsilon(\mathfrak{I})}$ is a seminorm over the interfacial bonds. Also note that
\begin{equation*}
  \|u'\|^2_{\ell^2_\varepsilon}
  =
  \|u'\|^2_{\ell^2_\varepsilon(\mathfrak{A})}
  +
  \|u'\|^2_{\ell^2_\varepsilon(\mathfrak{C})}
  +
  \|u'\|^2_{\ell^2_\varepsilon(\mathfrak{I})}.
\end{equation*}

\subsection{Variations}
The first variation, $\delta\E^{\text{QNL}}(y)[u]$, is
\begin{align*}
  \delta\E^{\text{QNL}}(y)[u]
  &=
  \varepsilon\sum_{\ell=1}^N\frac{\phi'(\|y'_\ell\|)}{\|y'_\ell\|}y'_\ell\cdot u'_\ell
  +
  \varepsilon\sum_{\ell=1}^K\frac{\phi'(\|y'_{\ell+1}+y'_\ell\|)}{\|y'_{\ell+1}+y'_\ell\|}(y'_{\ell+1}+y'_\ell)\cdot(u'_{\ell+1}+u'_\ell)\\
  &\quad+
  \varepsilon\sum_{\ell=K+2}^N2\,\frac{\phi'(2\|y'_\ell\|)}{\|y'_\ell\|}y'_\ell\cdot u'_\ell
  +
  \varepsilon\,\frac{\phi'(2\|y'_1\|)}{\|y'_1\|}y'_1\cdot u'_1
  +
  \varepsilon\,\frac{\phi'(2\|y'_{K+1}\|)}{\|y'_{K+1}\|}y'_{K+1}\cdot u'_{K+1},
\end{align*}
and the second variation is
\begin{align*}
  \delta^2\E^{\text{QNL}}(y)[u,v]
  &=
  \varepsilon\sum_{\ell=1}^N
  u'_\ell\cdot\left(\phi''(\|y'_\ell\|)P_\ell+\frac{\phi'(\|y'_\ell\|)}{\|y'_\ell\|}(I-P_\ell)\right)v'_\ell\\
  &\quad+
  \varepsilon\sum_{\ell=1}^K
  (u'_{\ell+1}+u'_\ell)\cdot\left(\phi''(\|y'_{\ell+1}+y'_\ell\|)\tilde{P}_\ell+\frac{\phi'(\|y'_{\ell+1}+y'_\ell\|)}{\|y'_{\ell+1}+y'_\ell\|}(I-\tilde{P}_\ell)\right)(v'_{\ell+1}+v'_\ell)\\
  &\quad+
  \varepsilon\sum_{\ell=K+2}^N
  u'_\ell\cdot\left(4\,\phi''(2\|y'_\ell\|)P_\ell+2\,\frac{\phi'(2\|y'_\ell\|)}{\|y'_\ell\|}(I-P_\ell)\right)v'_\ell\\
  &\quad+
  \frac{\varepsilon}{2}\,u'_1\cdot\left(4\,\phi''(2\|y'_1\|)P_1+2\,\frac{\phi'(2\|y'_1\|)}{\|y'_1\|}(I-P_1)\right)v'_1\\
  &\quad+
  \frac{\varepsilon}{2}\,u'_{K+1}\cdot\left(4\,\phi''(2\|y'_{K+1}\|)P_{K+1}+2\,\frac{\phi'(2\|y'_{K+1}\|)}{\|y'_{K+1}\|}(I-P_{K+1})\right)v'_{K+1}.
\end{align*}
Assuming again that the strains $F_1$ and $F_2$ defined in \eqref{eq:strains} are independent of $\ell$, we can rewrite the second variation as
\begin{align*}
  \delta^2\E^{\text{QNL}}(y)[u,v]
  &=
  \varepsilon\sum_{\ell=1}^N
  u'_\ell\cdot\left(\phi''(F_1)P_\ell+\frac{\phi'(F_1)}{F_1}(I-P_\ell)\right)v'_\ell\\
  &\quad+
  \varepsilon\sum_{\ell=1}^K
  (u'_{\ell+1}+u'_\ell)\cdot\left(\phi''(2F_2)\tilde{P}_\ell+\frac{\phi'(2F_2)}{2F_2}(I-\tilde{P}_\ell)\right)(v'_{\ell+1}+v'_\ell)\\
  &\quad+
  \varepsilon\sum_{\ell=K+2}^N
  u'_\ell\cdot\left(4\,\phi''(2F_1)P_\ell+2\,\frac{\phi'(2F_1)}{F_1}(I-P_\ell)\right)v'_\ell\\
  &\quad+
  \frac{\varepsilon}{2}\,u'_1\cdot\left(4\,\phi''(2F_1)P_1+2\,\frac{\phi'(2F_1)}{F_1}(I-P_1)\right)v'_1\\
  &\quad+
  \frac{\varepsilon}{2}\,u'_{K+1}\cdot\left(4\,\phi''(2F_1)P_{K+1}+2\,\frac{\phi'(2F_1)}{F_1}(I-P_{K+1})\right)v'_{K+1}.
\end{align*}
If we now apply the same manipulations preceding the derivation of \eqref{eq:E^a''}, rearrange the sums suitably, and collect interface terms, we obtain
\begin{align}
  \label{eq:E^QNL''}
  \delta^2\E^{\text{QNL}}(y)[u,v]
  &=
  \varepsilon\sum_{\ell=2}^K
  u'_\ell\cdot\left(\bigr(\phi''(F_1)+4\,\phi''(2F_1)\bigr)P_\ell+\frac{\phi'(F_1)}{F_1}(I_2-P_\ell)\right)v'_\ell\notag\\
  &\quad+
  \varepsilon\sum_{\ell=1}^K
  \biggr[
  (u'_{\ell+1}+u'_\ell)\cdot\left(\frac{\phi'(2F_1)}{2F_1}(I_2-\tilde{P}_\ell)\right)(v'_{\ell+1}+v'_\ell)
  -
  \varepsilon^2u''_{\ell+1}\cdot\left(\phi''(2F_1)\tilde{P}_\ell\right)v''_{\ell+1}
  \biggr]\notag\\
  &\quad+
  2\,\varepsilon\sum_{\ell=2}^K
  u'_\ell\cdot\left(\phi''(2F_1)(\tilde{P}_\ell+\tilde{P}_{\ell-1}-2P_\ell)\right)v'_\ell\notag\\
  &\quad+
  \varepsilon\sum_{\ell=1}^K
  \biggr[
  (u'_{\ell+1}+u'_\ell)\cdot\left((\phi''(2F_2)-\phi''(2F_1))\tilde{P}_\ell\right)(v'_{\ell+1}+v'_\ell)\\
  &\qquad\qquad+
  (u'_{\ell+1}+u'_\ell)\cdot\left(\left(\frac{\phi'(2F_2)}{2F_2}-\frac{\phi'(2F_1)}{2F_1}\right)(I_2-\tilde{P}_\ell)\right)(v'_{\ell+1}+v'_\ell)
  \biggr]\notag\\
  &+
  \varepsilon\sum_{\ell=K+2}^N
  u'_\ell\cdot\left(\bigr(\phi''(F_1)+4\,\phi''(2F_1)\bigr)P_\ell+\frac{\phi'(F_1)+2\phi'(2F_1)}{F_1}\left(I_2-P_\ell\right)\right)v'_{\ell}\notag\\
  &+
  \varepsilon\,u'_1\cdot\left((\phi''(F_1)+4\,\phi''(2F_1))P_1+\frac{\phi'(F_1)+\phi'(2F_1)}{F_1}(I_2-P_1)\right)v'_1\notag\\
  &+
  \varepsilon\,u'_{K+1}\cdot\left((\phi''(F_1)+4\,\phi''(2F_1))P_{K+1}+\frac{\phi'(F_1)+\phi'(2F_1)}{F_1}(I_2-P_{K+1})\right)v'_{K+1}\notag\\
  &+
  2\,\varepsilon\,u'_1\cdot\left(\phi''(2F_1)(\tilde{P}_1-P_1)\right)v'_1
  +
  2\,\varepsilon\,u'_{K+1}\cdot\left(\phi''(2F_1)(\tilde{P}_{K+1}-P_{K+1})\right)v'_{K+1}.\notag
\end{align}
Thus we see that the second variation can be roughly decomposed into the second variation of the atomistic energy (first five lines in \eqref{eq:E^QNL''}; cf.~\eqref{eq:E^a''}) and of the Cauchy--Born energy (the sixth line; cf.~\eqref{eq:E^CB''}), with some terms arising at the interfaces between the atomistic and continuum regions (last three lines).

\subsection{Stability of the quasi-nonlocal approximation}
Using the second variation of $\delta^2\E^{\text{QNL}}$ given in \eqref{eq:E^QNL''}, we can now give sharp stability results similar to those in Sections \ref{sec:stab_linear} and \ref{sec:stab_circ} for the atomistic and Cauchy--Born models.

For the 1-D chain, $y_F\in\tilde{\Y}$, we have that $P_\ell=\tilde{P}_\ell=P$ for all $\ell$ and $F_1=F_2=F$. In addition, for the constrained chain, for which the displacements satisfy $u\in\tilde{\U}$, we also have $(I-P)u_\ell=0$ for all $\ell$. Therefore, we immediately have the following stability result for the 1-D constrained chain. This result is consistent with that derived in~\cite{doblusort:qce.stab}.
\begin{theorem}
  \label{thm:stab_QNL_1D1D}
  Let $y_F\in\tilde{\Y}$ denote the 1-D configuration of atoms with nearest neighbor interatomic spacing $F\varepsilon$. If $\phi''(2F)\le0$, then
  \begin{equation*}
    \inf_{u\in\tilde{\U}\setminus\{0\}}
    \frac{\delta^2\E^{\text{QNL}}(y_F)[u,u]}{\|u'\|_{\ell^2_\varepsilon}^2}
    =
    \phi''(F)+4\,\phi''(2F).
  \end{equation*}
\end{theorem}
\begin{proof}
  Using expression \eqref{eq:E^QNL''} for the second variation of $\E^{\text{QNL}}$, we have
  \begin{align*}
    \delta^2\E^{\text{QNL}}(y_F)[u,u]
    &=
    \varepsilon\sum_{\ell=1}^N
    (\phi''(F)+4\,\phi''(2F))\|u'_\ell\|^2
    -
    \varepsilon\sum_{\ell=1}^K
    \varepsilon^2\phi''(2F)\|u''_{\ell+1}\|^2\notag\\
    &\ge
    (\phi''(F)+4\,\phi''(2F))\|u'\|^2_{\ell^2_\varepsilon}.
  \end{align*}
  The sharp estimate for the infimum follows by choosing a test function $u\in\tilde{\U}$ supported away from the atomistic region.
\end{proof}

For the 1-D unconstrained chain, we still have $P_\ell=\tilde{P}_\ell=P$ for all $\ell$ and $F_1=F_2=F$, but since the displacements are no longer one-dimensional, we no longer have $(I-P)u_\ell=0$ for all $\ell$. The stability result for the 1-D unconstrained chain follows.
\begin{theorem}
  \label{thm:stab_QNL_1D2D}
  Let $y_F\in\tilde{\Y}$ denote the 1-D configuration of atoms with nearest neighbor interatomic spacing $F\varepsilon$. If $\phi'(2F)\ge0$ and $\phi''(2F)\le0$, then
  \begin{align*}
    \delta^2\E^{\text{QNL}}(y_F)[u,u]
    &\ge
    \min{\left\{\phi''(F)+4\,\phi''(2F),\frac{\phi'(F)}{F}\right\}}\|u'\|^2_{\ell^2_\varepsilon(\mathfrak{A})}\notag\\
    &\quad+
    \min{\left\{\phi''(F)+4\,\phi''(2F),\frac{\phi'(F)+2\phi'(2F)}{F}\right\}}\|u'\|^2_{\ell^2_\varepsilon(\mathfrak{C})}\notag\\
    &\quad+
    \min{\left\{\phi''(F)+4\,\phi''(2F),\frac{\phi'(F)+\phi'(2F)}{F}\right\}}\|u'\|^2_{\ell^2_\varepsilon(\mathfrak{I})}\notag\\
    &\ge
    \min{\left\{\phi''(F)+4\,\phi''(2F),\frac{\phi'(F)}{F}\right\}}\|u'\|^2_{\ell^2_\varepsilon}
  \end{align*}
  for all $u\in\U$. In addition, as $\varepsilon\to0$,
  \begin{equation*}
    \inf_{u\in\U\setminus\{0\}}
    \frac{\delta^2\E^{\text{QNL}}(y_F)[u,u]}{\|u'\|_{\ell^2_\varepsilon}^2}
    =
    \min{\left\{\phi''(F)+4\,\phi''(2F),\frac{\phi'(F)}{F}+\mathcal{O}(\varepsilon)\right\}}.
  \end{equation*}
\end{theorem}
\begin{proof}
  Using expression \eqref{eq:E^QNL''} for the second variation of $\E^{\text{QNL}}$, we have
  \begin{align*}
    \delta^2\E^{\text{QNL}}(y_F)[u,u]
    &=
    \varepsilon\sum_{\ell=2}^K
    \biggr[(\phi''(F)+4\,\phi''(2F))\|P u'_\ell\|^2+\frac{\phi'(F)}{F}\|(I_2-P)u'_\ell\|^2\biggr]\\
    &\qquad+
    \varepsilon\sum_{\ell=1}^K
    \biggr[
    \frac{\phi'(2F)}{2F}\|(I_2-P)(u'_{\ell+1}+u'_\ell)\|^2
    -
    \varepsilon^2\phi''(2F)\|P u''_{\ell+1}\|^2
    \biggr]\\
    &\qquad+
    \varepsilon\sum_{\ell=K+2}^N
    \biggr[(\phi''(F)+4\,\phi''(2F)\bigr)\|P u'_\ell\|^2+\frac{\phi'(F)+2\,\phi'(2F)}{F}\|(I_2-P)u'_\ell\|^2\biggr]\\
    &\qquad+
    \varepsilon\biggr[(\phi''(F)+4\,\phi''(2F))\|P u'_1\|^2+\frac{\phi'(F)+2\,\phi'(2F)}{F}\|(I_2-P)u'_1\|^2\biggr]\\
    &\qquad+
    \varepsilon\biggr[(\phi''(F)+4\,\phi''(2F))\|P u'_{K+1}\|^2+\frac{\phi'(F)+2\,\phi'(2F)}{F}\|(I_2-P)u'_{K+1}\|^2\biggr].
  \end{align*}
  Due to the assumptions $\phi'(2F)\ge0$ and $\phi''(2F)\le0$, we can drop the second sum and estimate the remaining terms to get
  \begin{align*}
    \delta^2\E^{\text{QNL}}(y_F)[u,u]
    &\ge
    \min{\left\{\phi''(F)+4\,\phi''(2F),\frac{\phi'(F)}{F}\right\}}\|u'\|^2_{\ell^2_\varepsilon(\mathfrak{A})}\\
    &\quad+
    \min{\left\{\phi''(F)+4\,\phi''(2F),\frac{\phi'(F)+2\phi'(2F)}{F}\right\}}\|u'\|^2_{\ell^2_\varepsilon(\mathfrak{C})}\\
    &\quad+
    \min{\left\{\phi''(F)+4\,\phi''(2F),\frac{\phi'(F)+\phi'(2F)}{F}\right\}}\|u'\|^2_{\ell^2_\varepsilon(\mathfrak{I})}.
  \end{align*}
  Since $\phi'(2F)\ge0$, the smallest of the minima above is the first one, and the first part of the theorem follows.

  To obtain the sharp result for the infimum, consider again the displacements $\tilde{u}\in\tilde\U$ and $\hat{u}\in\U$ as in the proofs of Theorems \ref{thm:stab_CB_2D} and \ref{thm:stab_a_2D}. If the support of $\tilde{u}\in\tilde\U$ is chosen so that $\tilde{u}''_{\ell+1}=0$ for $\ell=1,\dots,K$, then we obtain
  \begin{equation*}
    \delta^2\E^{\text{QNL}}(y_F)[\tilde{u},\tilde{u}]
    =
    (\phi''(F)+4\,\phi''(2F))\|\tilde{u}'\|^2_{\ell^2_\varepsilon}.
  \end{equation*}
  On the other hand, if the zig-zag displacement $\hat{u}\in\U$ is chosen so that its support is contained in the atomistic region (atoms $2,\dots,K-1$ to be specific), then there is no contribution to $\delta^2\E^{\text{QNL}}$ from the continuum or interfacial regions, however, there is a $\mathcal{O}(\varepsilon)$ contribution from the $\hat{u}''_{\ell+1}$ and $\hat{u}'_{\ell+1}+\hat{u}'_\ell$ terms near the boundary of the support. We then have
  \begin{equation*}
    \delta^2\E^{\text{QNL}}(y_F)[\hat{u},\hat{u}]
    =
    \left(\frac{\phi'(F)}{F}+\mathcal{O}(\varepsilon)\right)\|\hat{u}'\|^2_{\ell^2_\varepsilon},
  \end{equation*}
  and the theorem follows.
\end{proof}

Finally, in the case of the uniform circular chain, none of the projection operators are the same, nor $F_1=F_2$. However, using the techniques used in Section \ref{sec:stab_circ}, we have the following result.
\begin{theorem}
  \label{thm:stab_QNL_2D}
  Let $y_F\in\Y$ denote the uniform circular configuration of $N$ atoms with nearest neighbor interatomic spacing $F\varepsilon$. If $\phi'(2F)\ge0$ and $\phi''(2F)\le0$, then
  \begin{align*}
    \delta^2\E^{\text{QNL}}(y_F)[u,u]
    &\ge
    \min{\left\{\phi''(F)+4\,\phi''(2F),\frac{\phi'(F)}{F}\right\}}\|u'\|^2_{\ell^2_\varepsilon(\mathfrak{A})}\notag\\
    &\quad+
    \min{\left\{\phi''(F)+4\,\phi''(2F),\frac{\phi'(F)+2\phi'(2F)}{F}\right\}}\|u'\|^2_{\ell^2_\varepsilon(\mathfrak{C})}\notag\\
    &\quad+
    \min{\left\{\phi''(F)+4\,\phi''(2F),\frac{\phi'(F)+\phi'(2F)}{F}\right\}}\|u'\|^2_{\ell^2_\varepsilon(\mathfrak{I})}\notag\\
    &\quad-
    \varepsilon^2(4\pi^2+4C_\phi)\|u'\|^2_{\ell^2_\varepsilon(\mathfrak{A})}\notag\\
    &\quad-
    \varepsilon\,(2\pi|\phi''(2F)|+2\,\varepsilon\,C_\phi)\|u'\|^2_{\ell^2_\varepsilon(\mathfrak{I})}\notag\\
    &\ge
    \min{\left\{\phi''(F)+4\,\phi''(2F),\frac{\phi'(F)}{F}\right\}}\|u'\|^2_{\ell^2_\varepsilon}\\
    &\quad-
    \varepsilon\,\max{\bigr\{2\pi|\phi''(2F)|+2\,\varepsilon\,C_\phi,4\,\varepsilon(\pi^2+C_\phi)\bigr\}}\|u'\|^2_{\ell^2_\varepsilon}\notag
  \end{align*}
  for all $u\in\U$, where $C_\phi$ is the Lipschitz constant defined in \eqref{eq:C_phi}. In addition, as $\varepsilon\to0$,
  \begin{equation*}
    \inf_{u\in\U\setminus\{0\}}
    \frac{\delta^2\E^{\text{QNL}}(y_F)[u,u]}{\|u'\|_{\ell^2_\varepsilon}^2}
    =
    \min{\left\{\phi''(F)+4\,\phi''(2F),\frac{\phi'(F)}{F}\right\}}
    +
    \mathcal{O}(\varepsilon).
  \end{equation*}
\end{theorem}
\begin{proof}
  Using expression \eqref{eq:E^QNL''} for the second variation of $\E^{\text{QNL}}$, we have
  \begin{align*}
    \delta^2\E^{\text{QNL}}(y_F)[u,u]
    &=
    \varepsilon\sum_{\ell=2}^K
    \biggr[(\phi''(F)+4\,\phi''(2F))\|P_\ell u'_\ell\|^2
    +
    \frac{\phi'(F)}{F}\|(I_2-P_\ell)u'_\ell\|^2\\
    &\qquad\qquad+
    2\,u'_\ell\cdot(\phi''(2F)(\tilde{P}_\ell+\tilde{P}_{\ell-1}-2P_\ell))u'_\ell\biggr]\\
    &\qquad+
    \varepsilon\sum_{\ell=1}^K
    \biggr[
    \frac{\phi'(2F)}{2F}\|(I_2-\tilde{P}_\ell)(u'_{\ell+1}+u'_\ell)\|^2
    -
    \varepsilon^2\phi''(2F)\|\tilde{P}_\ell u''_{\ell+1}\|^2
    \biggr]\\
    &\qquad+
    \varepsilon\sum_{\ell=1}^K
    \biggr[(\phi''(2F_2)-\phi''(2F))\|\tilde{P}_\ell(u'_{\ell+1}+u'_\ell)\|^2\\
    &\qquad\qquad
    \left(\frac{\phi'(2F_2)}{2F_2}-\frac{\phi'(2F)}{2F}\right)\|(I_2-\tilde{P}_\ell)(u'_{\ell+1}+u'_\ell)\|^2\biggr]\\
    &\qquad+
    \varepsilon\sum_{\ell=K+2}^N
    \biggr[(\phi''(F)+4\,\phi''(2F))\|P_\ell u'_\ell\|^2
    +
    \frac{\phi'(F)+2\phi'(2F)}{F}\|(I_2-P_\ell)u'_\ell\|^2\biggr]\\
    &\qquad+
    \varepsilon\biggr[(\phi''(F)+4\,\phi''(2F))\|P u'_1\|^2+\frac{\phi'(F)+2\,\phi'(2F)}{F}\|(I_2-P)u'_1\|^2\biggr]\\
    &\qquad+
    \varepsilon\biggr[(\phi''(F)+4\,\phi''(2F))\|P u'_{K+1}\|^2+\frac{\phi'(F)+2\,\phi'(2F)}{F}\|(I_2-P)u'_{K+1}\|^2\biggr]\\
    &\qquad+
    2\,\varepsilon\,u'_1\cdot\left(\phi''(2F)(\tilde{P}_1-P_1)\right)u'_1
    +
    2\,\varepsilon\,u'_{K+1}\cdot\left(\phi''(2F)(\tilde{P}_{K+1}-P_{K+1})\right)u'_{K+1}.
  \end{align*}
  Due to the assumptions $\phi'(2F)\ge0$ and $\phi''(2F)\le0$, we can drop the second sum and estimate the remaining terms to get
  \begin{align*}
    \delta^2\E^{\text{QNL}}(y_F)[u,u]
    &\ge
    \min{\left\{\phi''(F)+4\,\phi''(2F),\frac{\phi'(F)}{F}\right\}}\|u'\|^2_{\ell^2_\varepsilon(\mathfrak{A})}\\
    &\qquad+
    \min{\left\{\phi''(F)+4\,\phi''(2F),\frac{\phi'(F)+2\phi'(2F)}{F}\right\}}\|u'\|^2_{\ell^2_\varepsilon(\mathfrak{C})}\\
    &\qquad+
    \min{\left\{\phi''(F)+4\,\phi''(2F),\frac{\phi'(F)+\phi'(2F)}{F}\right\}}\|u'\|^2_{\ell^2_\varepsilon(\mathfrak{I})}\\
    &\qquad-
    4\,\varepsilon^2\pi^2\|u'\|^2_{\ell^2_\varepsilon(\mathfrak{A})}
    -
    4\,\varepsilon^2C_\phi\left(\|u'\|^2_{\ell^2_\varepsilon(\mathfrak{A})}+\|u'\|^2_{\ell^2_\varepsilon(\mathfrak{I})}\right)
    -
    2\,\varepsilon\pi|\phi''(2F)|\|u'\|^2_{\ell^2_\varepsilon(\mathfrak{I})},
  \end{align*}
  and the first part of the theorem follows.

  To obtain the sharp result for the infimum, consider again the displacements $\tilde{u}\in\U$ and $\hat{u}\in\U$ as in the proof of Theorems \ref{thm:stab_CB} and \ref{thm:stab_a}. Using the expansion displacement $\tilde{u}\in\U$, we obtain
  \begin{equation*}
    \delta^2\E^{\text{QNL}}(y_F)[\tilde{u},\tilde{u}]
    =
    \left(\phi''(F)+4\,\phi''(2F)+\mathcal{O}(\varepsilon)\right)\|\tilde{u}'\|^2_{\ell^2_\varepsilon},
  \end{equation*}
  where the $\mathcal{O}(\varepsilon)$ term arises from the interfacial terms. We note that even if the support of $\tilde{u}$ was restricted to a smaller set, a $\mathcal{O}(\varepsilon)$ term would still arise due to the terms from the boundary of the support.

  On the other hand, if the support of the zig-zag displacement $\hat{u}\in\U$ is again chosen so that there is no contribution to $\delta^2\E^{\text{QNL}}$ from the continuum or interfacial regions, we can obtain
  \begin{equation*}
    \delta^2\E^{\text{QNL}}(y_F)[\hat{u},\hat{u}]
    =
    \left(\frac{\phi'(F)}{F}+\mathcal{O}(\varepsilon)\right)\|\hat{u}'\|^2_{\ell^2_\varepsilon},
  \end{equation*}
  where again the $\mathcal{O}(\varepsilon)$ term arises from the term on the boundary of the support of $\hat{u}$. The result now follows.
\end{proof}

\subsection{Modeling error for the quasi-nonlocal approximation}
We next estimate the modeling error of the quasi-nonlocal approximation. Let $f\in\U$ denote an external force applied to a deformation $y_F\in\Y$ and let $u^{\text{a}}\in\U$ and $u^{\text{QNL}}\in\U$ solve the linearized equations
\begin{alignat*}
  \delta\E^{\text{a}}(y_F)[v]+\delta^2\E^{\text{a}}(y_F)[u^{\text{a}},v]
  &=
  \langle f,v\rangle
  &\qquad\text{ for all }v\in\U,\\
  \delta\E^{\text{QNL}}(y_F)[v]+\delta^2\E^{\text{QNL}}(y_F)[u^{\text{QNL}},v]
  &=
  \langle f,v\rangle
  &\qquad\text{ for all }v\in\U.
\end{alignat*}
The modeling error of the quasi-nonlocal approximation, $\tau$, is again given via the duality relationship
\begin{equation}
  \label{eq:trunc_def1}
  \langle\tau,v\rangle
  :=
  \left\{\delta\E^{\text{QNL}}(y_F)[v]-\delta\E^{\text{a}}(y_F)[v]\right\}
  +
  \left\{\delta^2\E^{\text{QNL}}(y_F)[u^{\text{a}},v]-\delta^2\E^{\text{a}}(y_F)[u^{\text{a}},v]\right\}
  \quad\text{ for all }v\in\U.
\end{equation}
Yet again, the solution $u^{\text{QNL}}$ will play no role in the analysis of the modeling error, so to simplify the notation in the rest of this section we will sometimes suppress the superscript and simply write $u$ instead of $u^{\text{a}}$. However, in the statements of the theorems the proper notation will be used.

We compute that we have for all $v\in\U$
\begin{align}
  \label{eq:ghost_QNL}
  \delta\E^{\text{QNL}}(y_F)[v]-\delta\E^{\text{a}}(y_F)[v]
  &=
  \varepsilon\sum_{\ell=K+2}^N\left[\left(\frac{\phi'(2F_1)}{2F_1}-\frac{\phi'(2F_2)}{2F_2}\right)\left(4y'_{F,\ell}
  +\varepsilon^2y'''_{F,\ell+1}\right)-\varepsilon^2\frac{\phi'(2F_1)}{2F_1}y'''_{F,\ell+1}\right]\cdot v'_\ell\notag\\
  &\qquad+
  \varepsilon\left[\left(\frac{\phi'(2F_1)}{2F_1}-\frac{\phi'(2F_2)}{2F_2}\right)\left(2y'_{F,1}-\varepsilon\,y''_{F,1}\right)
  +
  \varepsilon\frac{\phi'(2F_1)}{2F_1}y''_{F,1}\right]\cdot v'_1\\
  &\qquad+
  \varepsilon\left[\left(\frac{\phi'(2F_1)}{2F_1}-\frac{\phi'(2F_2)}{2F_2}\right)\left(2y'_{F,K+1}
  +
  \varepsilon\,y''_{F,K+2}\right)-\varepsilon\frac{\phi'(2F_1)}{2F_1}y''_{F,K+2}\right]\cdot v'_{K+1}\notag
\end{align}
and thus observe that
\begin{equation*}
  \delta\E^{\text{QNL}}(y_F)-\delta\E^{\text{a}}(y_F)\equiv0
  \qquad\text{if $y_F$ is a linear chain since $F_1=F_2$ and $y''_{F}\equiv y'''_{F}\equiv0$,}
\end{equation*}
but that
\begin{equation*}
  \delta\E^{\text{QNL}}(y_F)-\delta\E^{\text{a}}(y_F)\not\equiv0
  \qquad\text{if $y_F$ is a circular chain.}
\end{equation*}

After manipulations analogous to those in Section \ref{sec:truncation}, we can see that
\begin{align}
  \label{eq:trunc_def_qnl}
  \delta^2\E^{\text{QNL}}(y_F)[u,v]-\delta^2\E^{\text{a}}(y_F)[u,v]
  &=
  -\varepsilon\sum_{\ell=K+2}^N
  \biggr[\varepsilon^2\left(\tilde{A}_\ell u'''_{\ell+1}+\frac{\tilde{A}_\ell-\tilde{A}_{\ell-1}}{\varepsilon}u''_\ell\right)+2\left(\tilde{A}_\ell+\tilde{A}_{\ell-1}-2A_\ell\right)u'_\ell\biggr]\cdot v'_\ell\notag\\
  &\qquad+
  \varepsilon^2u''_1\cdot\tilde{A}_0v'_1
  -
  \varepsilon^2u''_{K+2}\cdot\tilde{A}_{K+1} v'_{K+1}\\
  &\qquad+
  2\,\varepsilon\,u'_1\cdot\left(A_1-\tilde{A}_0\right)v'_1
  +
  2\,\varepsilon\,u'_{K+1}\cdot\left(A_{K+1}-\tilde{A}_{K+1}\right)v'_{K+1},\notag
\end{align}
where the matrices $A_\ell$ and $\tilde{A}_\ell$ are defined in \eqref{eq:A_l}. We now, analogously to Theorems \ref{thm:truncation_1D1D} and \ref{thm:truncation_1D2D}, obtain the following results for 1-D constrained (Theorem \ref{thm:truncation_qnl_1D1D}) and 1-D unconstrained (Theorem \ref{thm:truncation_qnl_1D2D}) chains.
\begin{theorem}
  \label{thm:truncation_qnl_1D1D}
  Let $y_F\in\tilde{\Y}$ denote the 1-D configuration of atoms with nearest neighbor interatomic spacing $F\varepsilon$, let $f\in\tilde{\U}$, and let $u^{\text{a}}\in\tilde{\U}$ satisfy $\delta\E^{\text{a}}(y_F)[v]+\delta^2\E^{\text{a}}(y_F)[u^{\text{a}},v]=\langle f,v\rangle$ for all $v\in\tilde{\U}$. The modeling error of the quasi-nonlocal approximation, $\tau$, then satisfies the inequality
  \begin{equation*}
    \|\tau\|_*
    \le
    \varepsilon\,|\phi''(2F)|\left[\biggr(\varepsilon\|(u^\text{a})''_1\|^2+\varepsilon\|(u^\text{a})''_{K+2}\|^2\biggr)^{1/2}+\varepsilon\left(\varepsilon\sum_{\ell=K+2}^N\|(u^\text{a})'''_{\ell+1}\|^2\right)^{1/2}\right].
  \end{equation*}
\end{theorem}
\begin{proof}
  Since in this case $P_\ell=\tilde{P}_\ell=P$ for all $\ell$, $F_1=F_2=F$, and $(I_2-P)w=0$ for all $w\in\tilde{\U}$, the expression \eqref{eq:trunc_def_qnl} for the modeling error reduces to (recall that $\delta\E^{\text{QNL}}(y_F)-\delta\E^{\text{a}}(y_F)\equiv0$ for linear chains)
  \begin{align*}
    \langle\tau,v\rangle
    &=
    -\varepsilon\sum_{\ell=K+2}^N
    \varepsilon^2\phi''(2F)(u^\text{a})'''_{\ell+1}\cdot v'_\ell
    +
    \varepsilon^2\phi''(2F)\biggr((u^\text{a})''_1\cdot v'_1
    -
    (u^\text{a})''_{K+2}\cdot v'_{K+1}\biggr)\\
    &\le
    \varepsilon^2|\phi''(2F)|\left(\varepsilon\sum_{\ell=K+2}^N\|(u^\text{a})'''_{\ell+1}\|^2\right)^{1/2}\|v'\|_{\ell^2_\varepsilon(\mathfrak{C})}\\
    &\quad+
    \varepsilon|\phi''(2F)|\biggr(\varepsilon\|(u^\text{a})''_1\|^2+\varepsilon\|(u^\text{a})''_{K+2}\|^2\biggr)^{1/2}\|v'\|_{\ell^2_\varepsilon(\mathfrak{I})},
  \end{align*}
  and the result follows.
\end{proof}

\begin{theorem}
  \label{thm:truncation_qnl_1D2D}
  Let $y_F\in\tilde{\Y}$ denote the 1-D configuration of atoms with nearest neighbor interatomic spacing $F\varepsilon$, let $f\in\U$, and let $u^{\text{a}}\in\U$ satisfy $\delta\E^{\text{a}}(y_F)[v]+\delta^2\E^{\text{a}}(y_F)[u^{\text{a}},v]=\langle f,v\rangle$ for all $v\in\U$. The modeling error of the quasi-nonlocal approximation, $\tau$, then satisfies the inequality
  \begin{equation}
    \label{eq:trunc_qnl_1D2D}
    \|\tau\|_*
    \le
    \varepsilon\,\max{\left\{|\phi''(2F)|,\left|\frac{\phi'(2F)}{2F}\right|\right\}}\left[\biggr(\varepsilon\|(u^\text{a})''_1\|^2+\varepsilon\|(u^\text{a})''_{K+2}\|^2\biggr)^{1/2}+\varepsilon\left(\varepsilon\sum_{\ell=K+2}^N\|(u^\text{a})'''_{\ell+1}\|^2\right)^{1/2}\right].
  \end{equation}
\end{theorem}
\begin{proof}
  In this case, $P_\ell=\tilde{P}_\ell=P$ for all $\ell$ and $F_1=F_2=F$, but we do not, in general, have $(I_2-P)w=0$ for all $w\in\U$. Defining $A=\phi''(2F)P+\dfrac{\phi'(2F)}{2F}(I_2-P)$, the expression \eqref{eq:trunc_def_qnl} for the modeling error reduces to (recall that $\delta\E^{\text{QNL}}(y_F)-\delta\E^{\text{a}}(y_F)\equiv0$ for linear chains)
  \begin{align*}
    \langle\tau,v\rangle
    &=
    -\varepsilon\sum_{\ell=K+2}^N
    \varepsilon^2\left(A(u^\text{a})'''_{\ell+1}\right)\cdot v'_\ell
    +
    \varepsilon^2\biggr((u^\text{a})''_1\cdot Av'_1
    -
    (u^\text{a})''_{K+2}\cdot Av'_{K+1}\biggr)\\
    &\le
    \varepsilon^2\max{\left\{|\phi''(2F)|,\left|\frac{\phi'(2F)}{2F}\right|\right\}}\left(\varepsilon\sum_{\ell=K+2}^N\|(u^\text{a})'''_{\ell+1}\|^2\right)^{1/2}\|v'\|_{\ell^2_\varepsilon(\mathfrak{C})}\\
    &\quad+
    \varepsilon\,\max{\left\{|\phi''(2F)|,\left|\frac{\phi'(2F)}{2F}\right|\right\}}\biggr(\varepsilon\|(u^\text{a})''_1\|^2+\varepsilon\|(u^\text{a})''_{K+2}\|^2\biggr)^{1/2}\|v'\|_{\ell^2_\varepsilon(\mathfrak{I})},
  \end{align*}
  and the result follows.
\end{proof}

We now, analogously to Theorem \ref{thm:truncation_2D}, estimate the modeling error for the quasi-nonlocal approximation of a circular chain. In the proof, we will now need to estimate $\delta\E^{\text{QNL}}(y_F)-\delta\E^{\text{a}}(y_F)\not\equiv0$.
\begin{theorem}
  \label{thm:truncation_qnl_2D}
  Let $y_F\in\Y$ denote the uniform circular configuration of $N$ atoms with nearest neighbor interatomic spacing $F\varepsilon$, let $f\in\U$, and let $u^{\text{a}}\in\U$ satisfy $\delta\E^{\text{a}}(y_F)[v]+\delta^2\E^{\text{a}}(y_F)[u^{\text{a}},v]=\langle f,v\rangle$ for all $v\in\U$. If $\phi\in\C^3(0,\infty)$, then the modeling error of the quasi-nonlocal approximation, $\tau$, satisfies the inequality
  \begin{align*}
    \|\tau\|_*
    &\le
    (C_1\varepsilon^2+C_2\varepsilon^4)\left[\left(\varepsilon\sum_{\ell=K+2}^N\|(u^\text{a})'_\ell\|^2\right)^{1/2}+\left(\varepsilon\sum_{\ell=K+2}^N\|(u^\text{a})''_\ell\|^2\right)^{1/2}+\left(\varepsilon\sum_{\ell=K+2}^N\|(u^\text{a})'''_{\ell+1}\|^2\right)^{1/2}\right]\notag\\
    &\quad+
    (C_3\varepsilon+6\,C_\phi\varepsilon^2+C_\phi\varepsilon^3)
    \left[\biggr(\varepsilon\|(u^\text{a})'_1\|^2+\varepsilon\|(u^\text{a})'_{K+1}\|^2\biggr)^{1/2}
    +\biggr(\varepsilon\|(u^\text{a})''_1\|^2+\varepsilon\|(u^\text{a})''_{K+2}\|^2\biggr)^{1/2}\right]\\
    &\quad
    +
    C_\kappa\varepsilon^2
    +
    C_\mathfrak{I}\varepsilon^{3/2},
  \end{align*}
  where $C_1$ and $C_2$ are defined in \eqref{eq:C1C2}, $C_\phi$ is the Lipschitz constant defined in \eqref{eq:C_phi}, $C_\kappa$ is defined in \eqref{eq:C_kappa}, and
  \begin{gather}
    \label{eq:C3}
    C_3
    =
    \max{\left\{|\phi''(2F)|,\left|\frac{\phi'(2F)}{2F}\right|,2\pi\left|\phi''(2F)-\frac{\phi'(2F)}{2F}\right|\right\}},\notag\\
    C_\mathfrak{I}
    =
    \sqrt{2}\left(2\,\varepsilon\,C_\phi(1+\pi\varepsilon)F+\pi|\phi'(2F)|\right).
  \end{gather}
\end{theorem}
\begin{proof}
  Using \eqref{eq:ghost_QNL} and the fact that $\|y'_\ell\|=F$, $\|y''_\ell\|=2F\varepsilon^{-1}\sin(\pi\varepsilon)\le2F\pi$, and $\|y'''_\ell\|=4F\varepsilon^{-2}\sin^2(\pi\varepsilon)\le4F\pi^2$ for all $\ell$, we first estimate similarly as in \eqref{eq:ghost_estimate}
  \begin{equation*}
    \left|\delta\E^{\text{QNL}}(y_F)[v]-\delta\E^{\text{a}}(y_F)[v]\right|
    \le
    \varepsilon^2C_\kappa\|v'\|_{\ell^2_\varepsilon(\mathfrak{C})}
    +
    \varepsilon^{3/2}C_\mathfrak{I}\|v'\|_{\ell^2_\varepsilon(\mathfrak{I})}.
  \end{equation*}
  Using estimates similar to those in Theorem \ref{thm:truncation_2D}, we then have
  \begin{align*}
    \langle\tau,v\rangle
    &\le
    \varepsilon^2\left(\max{\left\{|\phi''(2F)|,\left|\frac{\phi'(2F)}{2F}\right|\right\}}+C_\phi\varepsilon^2\right)\left(\varepsilon\sum_{\ell=K+2}^N\|(u^\text{a})'''_{\ell+1}\|^2\right)^{1/2}\|v'\|_{\ell^2_\varepsilon(\mathfrak{C})}\\
    &\quad+
    \varepsilon^2\left(2\pi\left|\phi''(2F)-\frac{\phi'(2F)}{2F}\right|+4\pi\,C_\phi\varepsilon^2\right)\left(\varepsilon\sum_{\ell=K+2}^N\|(u^\text{a})''_\ell\|^2\right)^{1/2}\|v'\|_{\ell^2_\varepsilon(\mathfrak{C})}\\
    &\quad+
    \varepsilon^2\left(4\pi^2\left|\phi''(2F)-\frac{\phi'(2F)}{2F}\right|+12\,C_\phi\right)\left(\varepsilon\sum_{\ell=K+2}^N\|(u^\text{a})'_\ell\|^2\right)^{1/2}\|v'\|_{\ell^2_\varepsilon(\mathfrak{C})}\\
    &\quad+
    \varepsilon\left(\max{\left\{|\phi''(2F)|,\left|\frac{\phi'(2F)}{2F}\right|\right\}}+C_\phi\varepsilon^2\right)\biggr(\varepsilon\|(u^\text{a})''_1\|^2+\varepsilon\|(u^\text{a})''_{K+2}\|^2\biggr)^{1/2}\|v'\|_{\ell^2_\varepsilon(\mathfrak{I})}\\
    &\quad+
    \varepsilon\left(2\pi\left|\phi''(2F)-\frac{\phi'(2F)}{2F}\right|+6\,C_\phi\varepsilon\right)
    \biggr(\varepsilon\|(u^\text{a})'_1\|^2+\varepsilon\|(u^\text{a})'_{K+1}\|^2\biggr)^{1/2}
    \|v'\|_{\ell^2_\varepsilon(\mathfrak{I})}\\
    &\quad
    +
    \varepsilon^2C_\kappa\|v'\|_{\ell^2_\varepsilon(\mathfrak{C})}
    +
    \varepsilon^{3/2}C_\mathfrak{I}\|v'\|_{\ell^2_\varepsilon(\mathfrak{I})},
  \end{align*}
  and the result follows.
\end{proof}

\begin{remark}
  Notice that in all three scenarios, under the assumption that the number of interfaces between the atomistic and continuum regions remains constant, the modeling error is of order $\mathcal{O}(\varepsilon^{3/2})$ as $\varepsilon\to0$.
\end{remark}

\subsection{Error analysis for the quasi-nonlocal method}
As in Section \ref{sec:error}, we easily see that
\begin{equation}
  \label{eq:trunc_stab_qnl}
  \langle\tau,u^{\text{a}}-u^{\text{QNL}}\rangle
  =
  \delta^2\mathcal{E}^{\text{QNL}}(y_F)[u^{\text{a}}-u^{\text{QNL}},u^{\text{a}}-u^{\text{QNL}}],
\end{equation}
which can be combined with the stability results and the modeling errors from the previous sections to obtain estimates of the errors of the quasi-nonlocal method.

Let us first define the stability constants (see Theorems \ref{thm:stab_QNL_1D1D}, \ref{thm:stab_QNL_1D2D}, and \ref{thm:stab_QNL_2D})
\begin{equation*}
  \gamma_3
  =
  \phi''(F)+4\,\phi''(2F)
  \quad\text{ and }\quad
  \gamma_4
  =
  \min{\left\{\phi''(F)+4\,\phi''(2F),\,\dfrac{\phi'(F)}{F}\right\}}.
\end{equation*}
Combining \eqref{eq:trunc_stab_qnl} with Theorems \ref{thm:stab_QNL_1D1D} and \ref{thm:truncation_qnl_1D1D}, we then obtain the following theorem for the 1-D constrained chain.
\begin{theorem}
  \label{thm:error_qnl_1D1D}
  Let $y_F\in\tilde{\Y}$ denote the 1-D configuration of atoms with nearest neighbor interatomic spacing $F\varepsilon$. Given $f\in\tilde{\U}$, let $u^{\text{a}},u^{\text{QNL}}\in\tilde{\U}$ satisfy $\delta\E^{\text{a}}(y_F)[v]+\delta^2\E^{\text{a}}(y_F)[u^{\text{a}},v]=\langle f,v\rangle$ and $\delta\E^{\text{QNL}}(y_F)[v]+\delta^2\E^{\text{QNL}}(y_F)[u^{\text{QNL}},v]=\langle f,v\rangle$, respectively, for all $v\in\tilde{\U}$. If $\phi''(2F)\le0$ and $\gamma_3>0$, then
  \begin{equation*}
    \|(u^{\text{a}}-u^{\text{QNL}})'\|_{\ell^2_\varepsilon}
    \le
   \frac{\varepsilon|\phi''(2F)|}{\gamma_3}
   \left[\biggr(\varepsilon\|(u^\text{a})''_1\|^2+\varepsilon\|(u^\text{a})''_{K+2}\|^2\biggr)^{1/2}
   +\varepsilon\left(\varepsilon\sum_{\ell=K+2}^N\|(u^\text{a})'''_{\ell+1}\|^2\right)^{1/2}\right].
  \end{equation*}
\end{theorem}
Combining \eqref{eq:trunc_stab_qnl} with Theorems \ref{thm:stab_QNL_1D2D} and \ref{thm:truncation_qnl_1D2D}, we have the following theorem for the 1-D unconstrained chain.
\begin{theorem}
  \label{thm:error_qnl_1D2D}
  Let $y_F\in\tilde{\Y}$ denote the 1-D configuration of atoms with nearest neighbor interatomic spacing $F\varepsilon$. Given $f\in\U$, let $u^{\text{a}},u^{\text{QNL}}\in\U$ satisfy $\delta\E^{\text{a}}(y_F)[v]+\delta^2\E^{\text{a}}(y_F)[u^{\text{a}},v]=\langle f,v\rangle$ and $\delta\E^{\text{QNL}}(y_F)[v]+\delta^2\E^{\text{QNL}}(y_F)[u^{\text{QNL}},v]=\langle f,v\rangle$, respectively, for all $v\in\U$. If $\phi'(2F)\ge0$, $\phi''(2F)\le0$, and $\gamma_4>0$, then
  \begin{equation*}
   \|(u^{\text{a}}-u^{\text{QNL}})'\|_{\ell^2_\varepsilon}
   \le
  C\varepsilon
   \left[\biggr(\varepsilon\|(u^\text{a})''_1\|^2+\varepsilon\|(u^\text{a})''_{K+2}\|^2\biggr)^{1/2}
   +\varepsilon\left(\varepsilon\sum_{\ell=K+2}^N\|(u^\text{a})'''_{\ell+1}\|^2\right)^{1/2}\right],
  \end{equation*}
  where
  \begin{equation*}
    C
    =
    \gamma_4^{-1}\max{\left\{|\phi''(2F)|,\left|\frac{\phi'(2F)}{2F}\right|\right\}}.
  \end{equation*}
\end{theorem}
Finally, note that the stability ``constant'' in Theorem \ref{thm:stab_QNL_2D} depends on $\varepsilon$. To simplify the notation, we define
\begin{equation*}
  \gamma_\varepsilon
  =
  \gamma_4-\varepsilon\,\max{\bigr\{2\pi|\phi''(2F)|+2\,\varepsilon\,C_\phi,4\,\varepsilon(\pi^2+C_\phi)\bigr\}},
\end{equation*}
where $C_\phi$ is the Lipschitz constant defined in \eqref{eq:C_phi}. Combining \eqref{eq:trunc_stab_qnl} with Theorems \ref{thm:stab_QNL_2D} and \ref{thm:truncation_qnl_2D}, we then have the following theorem for the uniform circular chain.
\begin{theorem}
  \label{thm:error_qnl_2D}
  Let $y_F\in\Y$ denote the uniform circular configuration of $N$ atoms with nearest neighbor interatomic spacing $F\varepsilon$. Given $f\in\U$, let $u^{\text{a}},u^{\text{QNL}}\in\U$ satisfy $\delta\E^{\text{a}}(y_F)[v]+\delta^2\E^{\text{a}}(y_F)[u^{\text{a}},v]=\langle f,v\rangle$ and $\delta\E^{\text{QNL}}(y_F)[v]+\delta^2\E^{\text{QNL}}(y_F)[u^{\text{QNL}},v]=\langle f,v\rangle$, respectively, for all $v\in\U$. If $\phi\in\mathcal{C}^3(0,\infty)$, $\phi'(2F)\ge0$, $\phi''(2F)\le0$, and if $\gamma_4>0$ and $\varepsilon$ is small enough so that $\gamma_\varepsilon>0$, then
  \begin{align*}
    \|(u^{\text{a}}-&u^{\text{QNL}})'\|_{\ell^2_\varepsilon}
    \le\\
    &\gamma_\varepsilon^{-1}(C_1\varepsilon^2+C_2\varepsilon^4)\left[\left(\varepsilon\sum_{\ell=K+2}^N\|(u^\text{a})'_\ell\|^2\right)^{1/2}+\left(\varepsilon\sum_{\ell=K+2}^N\|(u^\text{a})''_\ell\|^2\right)^{1/2}+\left(\varepsilon\sum_{\ell=K+2}^N\|(u^\text{a})'''_{\ell+1}\|^2\right)^{1/2}\right]\\
    &+
    \gamma_\varepsilon^{-1}(C_3\varepsilon+6\,C_\phi\varepsilon^2
    +C_\phi\varepsilon^3)\left[\biggr(\varepsilon\|(u^\text{a})'_1\|^2
    +\varepsilon\|(u^\text{a})'_{K+1}\|^2\biggr)^{1/2}+\biggr(\varepsilon\|(u^\text{a})''_1\|^2
    +\varepsilon\|(u^\text{a})''_{K+2}\|^2\biggr)^{1/2}\right]\\
    &+
    \gamma_\varepsilon^{-1}(C_\kappa\varepsilon^2+C_\mathfrak{I}\varepsilon^{3/2}),
  \end{align*}
  where $C_1$ and $C_2$ are defined in \eqref{eq:C1C2}, $C_3$ and $C_\mathfrak{I}$ in \eqref{eq:C3}, $C_\phi$ in \eqref{eq:C_phi}, and $C_\kappa$ in \eqref{eq:C_kappa}.
\end{theorem}

\begin{remark}
  Stability estimates for quasi-nonlocal approximations of $\E^{\text{a,b}}(y)$ that include the bond-angle energy, $\E^{\text{QNL,b}}(y)$, can be obtained by including the increased stability~\eqref{increase}. This increased stability for the atomistic and Cauchy--Born model was analyzed in Section~\ref{sec:regularization}. The modeling error analysis for $\E^{\text{QNL,b}}(y)$ is identical to the modeling error analysis for $\E^{\text{QNL}}(y)$ since the bond-angle energy is not approximated. Corresponding error estimates for $\E^{\text{QNL,b}}(y)$ can then be obtained from these stability and modeling error estimates and are summarized in the conclusion section below.
\end{remark}

\section{Conclusion}
\label{sec:conclusion}
We summarize the results in this paper in the following tables for the case $\phi'(2F)\ge0$ and $\phi''(2F)\le0$. We first give results for the stability of the fully atomistic model. These stability bounds (and those that follow for the Cauchy--Born and quasi-nonlocal approximations) neglect $\mathcal{O}(\varepsilon)$ terms. We have modeled the strength of the bond-angle energy \eqref{eq:E^b} by $\alpha\ge0$.

\vspace{.3in}

\begin{center}
\begin{tabular}{|c|c|}
  \hline
  {\bf Atomistic model} & {\bf Stability}\\ \hline\hline
  1D (constrained)      & $\phi''(F)+4\,\phi''(2F)>0$\\ \hline
  1D (unconstrained)    & $\min{\left\{\phi''(F)+4\,\phi''(2F),\,\dfrac{\phi'(F)}{F}+ \dfrac{2\alpha}{F^2}\right\}}>0$\\ \hline
  circle                & $\min{\left\{\phi''(F)+4\,\phi''(2F),\,\dfrac{\phi'(F)}{F}+\dfrac{2\alpha\cos{\beta_\varepsilon}}{F^2}\right\}}>0$\\ \hline
\end{tabular}
\end{center}
\vspace{.3in}

We next give stabilty and error estimates for the Cauchy--Born approximation. The second-order error estimates $\mathcal{O}(\varepsilon^{2})$ for the Cauchy--Born approximation require that $\|(u^{\text{a}})'''\|_{\ell_\varepsilon^2}$ be bounded uniformly in $\varepsilon$, which is not the case for the approximation of atomistic configurations with defects.
\vspace{.3in}

\begin{center}
\begin{tabular}{|c|c|c|}
  \hline
  {\bf Cauchy--Born model} & {\bf Stability}                                                                                                                  & {\bf Error}\\ \hline\hline
  1D (constrained)         & $\phi''(F)+4\,\phi''(2F)>0$                                                                                                      & $\mathcal{O}(\varepsilon^2)$\\ \hline
  1D (unconstrained)       & $\min{\left\{\phi''(F)+4\,\phi''(2F),\,\dfrac{\phi'(F)+2\,\phi'(2F)}{F}+ \dfrac{2\alpha}{F^2}\right\}}>0$                        &  $\mathcal{O}(\varepsilon^2)$\\ \hline
  circle                   & $\min{\left\{\phi''(F)+4\,\phi''(2F),\,\dfrac{\phi'(F)+2\,\phi'(2F)}{F}+\dfrac{2\alpha\cos{\beta_\varepsilon}}{F^2}\right\}}>0$  &  $\mathcal{O}(\varepsilon^2)$\\ \hline
\end{tabular}
\end{center}
\vspace{.3in}

The lack of accuracy of the Cauchy--Born approximation for problems with defects is the motivation for the development of atomistic-to-continuum methods such as the quasi-nonlocal method which attain $\mathcal{O}(\varepsilon^{3/2})$ accuracy for problems with defects~\cite{dobs-qcf2,brian10,ortner:qnl1d}. We summarize below our stability and error estimates for the quasi-nonlocal approximation.
\vspace{.3in}

\begin{center}
\begin{tabular}{|c|c|c|}
  \hline
  {\bf QNL model}    & {\bf Stability}                                                                                                    & {\bf Error}\\ \hline\hline
  1D (constrained)   & $\phi''(F)+4\,\phi''(2F)>0$                                                                                        & $\mathcal{O}(\varepsilon^{3/2})$\\ \hline
  1D (unconstrained) & $\min{\left\{\phi''(F)+4\,\phi''(2F),\,\dfrac{\phi'(F)}{F}+\dfrac{2\alpha}{F^2}\right\}}>0$                        & $\mathcal{O}(\varepsilon^{3/2})$\\ \hline
  circle             & $\min{\left\{\phi''(F)+4\,\phi''(2F),\,\dfrac{\phi'(F)}{F}+\dfrac{2\alpha\cos{\beta_\varepsilon}}{F^2}\right\}}>0$ & $\mathcal{O}(\varepsilon^{3/2})$\\ \hline
\end{tabular}
\end{center}

\end{document}